\documentclass[11pt]{article}
\usepackage{float} 
\usepackage{times}
\usepackage{epsfig}
\usepackage{graphicx}
\usepackage{amsmath}
\usepackage{amssymb}
\usepackage{bm}
\usepackage{algorithm}
\usepackage{algorithmicx}
\usepackage{algpseudocode}
\usepackage{comment}
\usepackage{fullpage}
\usepackage{tikz}
\usepackage{pifont}
\usepackage{amsthm}
\usepackage{subfigure}
\usepackage{authblk}
\usepackage{appendix}

\allowdisplaybreaks

\newtheorem{lemma}{Lemma} 
\newtheorem{definition}{Definition} 

\newtheorem{theorem}{Theorem}

\newtheorem*{remark*}{Remark}

\begin{document}

\title{SketchyCoreSVD: SketchySVD from Random Subsampling of the Data Matrix}
\author{Chandrajit Bajaj, Yi Wang, Tianming Wang}
\affil{Department of Computer Science, Oden Institute for Computational Engineering and Sciences, University of Texas, Austin, Texas, USA.}
\date{\today}

\maketitle

\begin{abstract}
    We present a method called SketchyCoreSVD to compute the near-optimal rank $r$ SVD of a data matrix by building random sketches only from its subsampled columns and rows. We provide theoretical guarantees under incoherence assumptions, and validate the performance of our SketchyCoreSVD method on various large static and time-varying datasets. 
\end{abstract}

\section{Introduction}

Data matrices in practice often have fast decaying spectrum. Finding low rank approximation to a data matrix is a fundamental task in numerous applications. Due to the massive size of data in modern applications, it may not be possible to store the full data matrix, causing problems for classic solvers.

It turns out that randomized sketches of the data matrix, built from dimension reduction maps, suffice for low-rank approximation. The first one-pass algorithm based on random sketches appear in \cite{woolfe2008fast}. Such ideas is further advocated in the survey article \cite{halko2011finding}. Ever since, quite a few algorithms \cite{woodruff2014sketching,cohen2015dimensionality,boutsidis2016optimal,tropp2017practical} based on random sketches have been developed. \footnote{See Appendix \ref{sec:prior_algorithms} for a tabular summary description of prior sketch building algorithms.}

The most recent paper called SketchySVD \cite{tropp2019streaming} is shown to provide practical and consistent approximations compared to the predecessors. For a data matrix $\bm{A}$, SketchySVD finds the basis for the column/row space by QR decomposition of the left and right sketches while the first $k$ singular values are well preserved in the core sketch. Theoretical guarantees have been derived for sketches built from Gaussian maps. The more practical dimension reduction maps such as Scrambled Subsampled Randomized Fourier Transform (SSRFT) \cite{woolfe2008fast} and sparse sign matrices \cite{achlioptas2003database,li2006very} are shown to exhibit similar performances in practice. 

Building upon SketchySVD, we would like to claim that there is redundancy in the data matrix that can be further exploited when constructing sketches, thereby reducing the computation cost without sacrificing approximation accuracy. Such ``redundancy'' is characterized by incoherence\footnote{See Definition \ref{def:incoherence} in Section \ref{sec:alg_theory}.}, which is widely used in compressed sensing and matrix recovery. Several papers \cite{talwalkar2014matrix,gittens2011spectral,chiu2013sublinear} also establish guarantees for Nystr\"{o}m method \cite{williams2001using} based on the incoherence assumption. Overall, our contributions are the following:
\begin{itemize}
    \item Propose SketchyCoreSVD, a method to compute SketchySVD from subsampled columns and rows of the data matrix;
    \item Experiments verify that SketchyCoreSVD is able to reduce the computation cost of SketchySVD without sacrificing approximation accuracy;
    \item Theoretical guarantees based on incoherence assumption of the data matrix.
\end{itemize} 
\section{Algorithm}

\subsection{SketchyCoreSVD} 
Suppose the data matrix $\bm{A}\in\mathbb{R}^{M\times N}$ is huge in size, and its spectrum is decaying fast enough for compression. Let $r$ be a suitable choice of rank. The main goal is to find a near-optimal rank $r$ approximation to $\bm{A}$, via dimension reduction maps applied on subsampled columns and rows of $\bm{A}$. 
Choose sketch sizes $k$ and $s$ such that $r\leq k \leq s\leq \min\{m,n,m^{\prime},n^{\prime}\}$, where $m$, $n$, $m^{\prime}$ and $n^{\prime}$ satisfy
$$
\frac{m}{M}=\frac{n}{N}=p< 1,
\quad
\frac{m^{\prime}}{M}=\frac{n^{\prime}}{N}=q< 1,
\quad\text{and}\quad
p\leq q.
$$
Inspired by SketchySVD \cite{tropp2019streaming}, the steps of our method, termed SketchyCoreSVD, are the following.
\begin{enumerate}
\item Building sketches.
    \begin{itemize}
    \item Uniformly sample (without replacement) $m$ rows of $\bm{A}$. Denote the indices of the sampled rows to be $\Delta$. The map $\bm{\Gamma}\in\mathbb{R}^{k\times m}$ is applied to $\bm{A}^{(\Delta,:)}$ from the left,
    $$\bm{X}=\bm{\Gamma}\bm{A}^{(\Delta,:)}\in\mathbb{R}^{k\times N};$$
    \item Uniformly sample (without replacement) $n$ columns of $\bm{A}$. Denote the indices of the sampled columns to be $\Theta$. The map $\bm{\Omega}\in\mathbb{R}^{k\times n}$ is applied to $\bm{A}^{(:,\Theta)}$ from the right, 
    $$\bm{Y}=\bm{A}^{(:,\Theta)}\bm{\Omega}^*\in\mathbb{R}^{M\times k};$$
    \item Uniformly sample (without replacement) $m^{\prime}$ row indices and $n^{\prime}$ column indices, denoted by $\Delta^{\prime}$ and $\Theta^{\prime}$, respectively. Apply random matrix maps $\bm{\Phi}\in\mathbb{R}^{s\times m^{\prime}}$ and $\bm{\Psi}\in\mathbb{R}^{s\times n^{\prime}}$ to the intersection of $\bm{A}^{(\Delta^{\prime},\Theta^{\prime})}\in\mathbb{R}^{m^{\prime}\times n^{\prime}}$, i.e.,
    $$
    \bm{Z}=\bm{\Phi}\bm{A}^{(\Delta^{\prime},\Theta^{\prime})}\bm{\Psi}^*\in\mathbb{R}^{s\times s}.
    $$
    \end{itemize}
\item Computations.
\begin{itemize}
    \item We compute the QR decomposition of $\bm{X}^*$, 
    $$\bm{X}^*=\bm{P}\bm{R}_1,$$
    where $\bm{P}\in\mathbb{R}^{N\times k}$, and $\bm{R}_1\in\mathbb{R}^{k\times k}$;
    \item We compute the QR decomposition of $\bm{Y}$, 
    $$\bm{Y}=\bm{Q}\bm{R}_2,$$
    where $\bm{Q}\in\mathbb{R}^{M\times k}$, and $\bm{R}_2\in\mathbb{R}^{k\times k}$;
    \item We compute the core approximation,
    $$
    \bm{C} = (\bm{\Phi}\bm{Q}^{(\Delta^{\prime},:)})^{\dagger}\cdot\bm{Z}\cdot((\bm{\Psi}\bm{P}^{(\Theta^{\prime},:)})^{\dagger})^*\in\mathbb{R}^{k\times k};
    $$
    \item Denote $\hat{\bm{A}}=\bm{Q}\bm{C}\bm{P}^*$ as the initial approximation. The final near-optimal rank $r$ approximation to $\bm{A}$, denoted by $[[\hat{\bm{A}}]]_r$, is computed by
    $$
    \bm{Q}[[\bm{C}]]_r\bm{P}^*,
    $$
    where $[[\bm{C}]]_r$ is the best rank $r$ approximation to $\bm{C}$.
\end{itemize}
\end{enumerate}

While the storage cost is the same, SketchyCoreSVD has reduced computational cost. The following table shows the flops count of the two methods for Gaussian maps. If the same $k$ and $s$ are used, SketchyCoreSVD has lower complexity. The cost is reduced mainly when constructing the sketches. 

\begin{table}[H]
\begin{center}
\begin{tabular}{ccc}
\hline
\multicolumn{1}{|c|}{} & \multicolumn{1}{c|}{\text{SketchySVD}} & \multicolumn{1}{c|}{\text{SketchyCoreSVD}} \\
\hline
\multicolumn{1}{|c|}{$\bm{X}$} & \multicolumn{1}{c|}{$O(kMN)$} & \multicolumn{1}{c|}{$O(kpMN)$} \\
\hline
\multicolumn{1}{|c|}{$\bm{Y}$} & \multicolumn{1}{c|}{$O(kMN)$} & \multicolumn{1}{c|}{$O(kMN)$} \\
\hline
\multicolumn{1}{|c|}{$\bm{Z}$} & \multicolumn{1}{c|}{$O(sMN+s^2\min\{M,N\})$} & \multicolumn{1}{c|}{$O(sq^2MN+s^2q\min\{M,N\})$} \\
\hline
\multicolumn{1}{|c|}{QR of $\bm{X}$} & \multicolumn{1}{c|}{$O(k^2N)$} & \multicolumn{1}{c|}{$O(k^2N)$} \\
\hline
\multicolumn{1}{|c|}{QR of $\bm{Y}$} & \multicolumn{1}{c|}{$O(k^2M)$} & \multicolumn{1}{c|}{$O(k^2M)$} \\
\hline
\multicolumn{1}{|c|}{$\bm{C}$} & \multicolumn{1}{c|}{$O(k^2s+ks^2+ks(M+N))$} & \multicolumn{1}{c|}{$O(k^2s+ks^2+ksq(M+N))$} \\
\hline
\multicolumn{1}{|c|}{$[[\bm{C}]]_r$} & \multicolumn{1}{c|}{$O(k^3)$} & \multicolumn{1}{c|}{$O(k^3)$} \\
\hline
\multicolumn{1}{|c|}{$[[\hat{\bm{A}}]]_r$} & \multicolumn{1}{c|}{$O(kr(M+N)+r^2\min\{M,N\}+rMN)$} & \multicolumn{1}{c|}{$O(kr(M+N)+r^2\min\{M,N\}+rMN)$} \\
\hline
\end{tabular}
\end{center}
\end{table}

\subsection{Theoretical Guarantees}
\label{sec:alg_theory}

Our proofs generally follow the outline of \cite{tropp2019streaming}, though with substantial differences when dealing with subsamples of the data matrix. The fist step is to prove that $\bm{Q}$ and $\bm{P}$ capture the range and co-range of $\bm{A}$, which are expressed as
$$
\bm{A}\approx \bm{Q}\bm{Q}^*\bm{A}, \quad
\bm{A}\approx \bm{A}\bm{P}\bm{P}^*.
$$
These are proved in \cite{halko2011finding} for $\bm{Q}$ and $\bm{P}$ computed from the sketches of the full matrix. For our case, only a randomly selected subsets of the columns/rows are used. Thus we must impose certain conditions on the columns/rows of $\bm{A}$. 

The intuition is that the columns/rows are ``more or less the same''. Such a matrix property can be characterized by incoherence.  Suppose $\bm{A}\in\mathbb{R}^{M\times N}$ has an SVD of the following form
$$
\bm{U}\bm{\Sigma}\bm{V}^* = \left[\begin{array}{cc}
\bm{U}_1 & \bm{U}_2
\end{array}\right]
\left[\begin{array}{cc}
\bm{\Sigma}_1 & \\
              & \bm{\Sigma}_2
\end{array}\right]
\left[\begin{array}{cc}
\bm{V}_1 & \bm{V}_2
\end{array}\right]^*,
$$
where $[[\bm{A}]]_{r}=\bm{U}_1\bm{\Sigma}_1\bm{V}_1^*$ is the best rank $r$ approximation of $\bm{A}$.

\begin{definition}\label{def:incoherence}
$[[\bm{A}]]_{r}\in\mathbb{R}^{M\times N}$ is $(\mu,\nu)$-incoherent if 
$$
\max_i\left\{\|\bm{U}_1^{(i,:)}\|_2\right\}\leq\sqrt{\frac{\mu r}{M}}
\quad\text{and}\quad
\max_j\left\{\|\bm{V}_1^{(j,:)}\|_2\right\}\leq\sqrt{\frac{\nu r}{N}}.
$$
\end{definition}

Based on the incoherence assumption, we provide the following error bound.

\begin{theorem}\label{thm:range}
Suppose $\bm{A}=\bm{U}_1\bm{\Sigma}_1\bm{V}_1^*+\bm{U}_2\bm{\Sigma}_2\bm{V}_2^*$, where $[[\bm{A}]]_{r}=\bm{U}_1\bm{\Sigma}_1\bm{V}_1^*$, the best rank $r$ approximation of $\bm{A}$, is $(\mu,\nu)$-incoherent. Then for our SketchyCoreSVD algorithm with $k\geq r+4$, $m\geq 8\mu r\log r$ and $n\geq 8\nu r\log r$,
\begin{align*}
\max\left\{\|\bm{A}-\bm{Q}\bm{Q}^*\bm{A}\|_F,\|\bm{A}-\bm{A}\bm{P}\bm{P}^*\|_F\right\}
\leq\left(C_1(p,k,r)+1\right)\cdot\|\bm{\Sigma}_2\|_F+C_2(p,k,r)\cdot\|\bm{\Sigma}_2\|_2
\end{align*}
with probability at least $1-\frac{4}{r^3}-\frac{4}{k^3}$, where 
$$
C_1(p,k,r)=\sqrt{\frac{6e^2}{p}}\cdot\frac{k}{k-r+1}\cdot k^{\frac{3}{k-r+1}}
\quad\text{and}\quad
C_2(p,k,r)=\sqrt{\frac{36e^2}{p}}\cdot\frac{\sqrt{k\log k}}{k-r+1}\cdot k^{\frac{3}{k-r+1}}.
$$
\end{theorem}
\begin{remark*}
Note that $C_1(p,k,r)$ and $C_2(p,k,r)$ both decrease as $p$ or $k$ increases. Thus we are advised to use a not too small sampling ratio $p$, and a bigger sketch size $k$ whenever possible. 
\end{remark*}

Let $\bm{Q}$ and $\bm{P}$ be the basis computed by SketchyCoreSVD for the columns and rows, respectively. There exist $\mu^{\prime}\in[1,M]$ and $\nu^{\prime}\in[1,N]$ such that 
$$
\max_i\left\{\|\bm{Q}^{(i,:)}\|_2\right\}\leq\sqrt{\frac{\mu^{\prime}k}{M}}
\quad\text{and}\quad
\max_j\left\{\|\bm{P}^{(j,:)}\|_2\right\}\leq\sqrt{\frac{\nu^{\prime}k}{N}}.
$$
The existence of $\mu^{\prime}$ and $\nu^{\prime}$ can be easily shown.\footnote{For example, consider $\mu^{\prime}$. Since $\max_i\{\|\bm{Q}^{(i,:)}\|_2^2\}\leq k$, it follows that $\mu^{\prime}\leq M$. If $\mu^{\prime}<1$, then $\|\bm{Q}\|_F^2< M\cdot\frac{k}{M}=k$, contradicting to the fact that $\|\bm{Q}\|_F^2=k$.} We derive the following guarantee for the rank $r$ approximation $[[\hat{\bm{A}}]]_r$ computed by SketchyCoreSVD, provided that $m^{\prime}$ and $n^{\prime}$ are greater than some multiples of $\mu^{\prime}k\log k$ and $\nu^{\prime}k\log k$, respectively. 

\begin{theorem}
Condition on the success of Theorem~\ref{thm:range}, for our SketchyCoreSVD algorithm with $s\geq k+4$, $m^{\prime}\geq 8\mu^{\prime}k\log k$, and $n^{\prime}\geq 8\nu^{\prime}k\log k$, 
the initial approximation $\hat{\bm{A}}$ satisfies
$$    
\|\bm{A}-\hat{\bm{A}}\|_F\leq C_3(p,q,s,k,r)\cdot\|\bm{\Sigma}_2\|_F+C_4(p,q,s,k,r)\cdot\|\bm{\Sigma}_2\|_2
$$
with probability at least $1-\frac{4}{k^3}-\frac{6}{s^3}$, where
$$
C_3(p,q,s,k,r)=C_1(p,k,r)\cdot(\sqrt{3}C(q,s,k)+\sqrt{2}),
\quad
C_4(p,q,s,k,r)=C_2(p,k,r)\cdot(\sqrt{3}C(q,s,k)+\sqrt{2}), 
$$
$$
C(q,s,k) = \frac{6e^2}{q}\cdot\frac{s^{1+6/(s-k+1)}}{(s-k+1)^2}\cdot\left(\sqrt{s}+\sqrt{6\log s}\right)^2.
$$
Moreover, the final rank $r$ approximation satisfies
$$
\|\bm{A}-[[\hat{\bm{A}}]]_r\|_F\leq(2C_3(p,q,s,k,r)+1)\cdot\|\bm{\Sigma}_2\|_F+(2C_4(p,q,s,k,r))\cdot\|\bm{\Sigma}_2\|_2.
$$
\end{theorem}
\begin{remark*}
Note that $C(q,s,k)$ decreases as $q$ or $s$ increases. Thus we are advised to use a not too small sampling ratio $q$, and a bigger sketch size $s$ whenever possible. 

We admit the dependence on $\mu^{\prime}$ and $\nu^{\prime}$ is the deficiency of our current theoretical guarantees. Ideally, one should derive bounds for $\mu^{\prime}$ and $\nu^{\prime}$ based on $\mu$ and $\nu$. One way to get around this is to slightly modify the algorithm: after building the left and right sketches, we calculate the QR factorizations, and sample $m^{\prime}$ row indices and $n^{\prime}$ column indices based on the actual incoherence parameters, estimated from the row norms, of $\bm{Q}$ and $\bm{P}$, respectively. 

Nevertheless, we observe that $\mu^{\prime}=O(\mu)$, $\nu^{\prime}=O(\nu)$ in practice.\footnote{See Table~6 for some empirical evidences.} In the experiments shown in Section~\ref{sec:numerics}, we simply choose $m^{\prime}=m$ and $n^{\prime}=n$, i.e., $q=p$.
\end{remark*} 
\section{Numerical Experiments}
\label{sec:numerics}
If full SVD can be computed, the guide to choose the rank is to find the ``knee'' in the scree plot. At rank $r$,
$$
scree(r) = \frac{1}{\|\bm{A}\|_F^2}\sum_{i=r+1}\sigma_i^2(\bm{A}).
$$
For a rank $r$ approximation $[[\hat{\bm{A}}]]_r$ to $\bm{A}$, we measure its approximation error as
$$
err = \frac{\|[[\hat{\bm{A}}]]_r-\bm{A}\|_F^2}{\|\bm{A}\|_F^2}.
$$
We compare our SketchyCoreSVD method to SketchySVD \cite{tropp2019streaming}, which suggests one to choose $s\geq 2k+1$ and $k=O(r)$. The dimension reduction maps used are Gaussian maps. For our SketchyCoreSVD, we choose the sampling ratio $p$ such that $p\cdot\min\{M,N\}\geq s$, and $q=p$. 

The experiments are executed from MATLAB R2018a on a 64-bit Linux machine with 8 Intel i7-7700 CPUs at 3.60 GHz and 32 GB of RAM. The reported errors are averaged over 20 trials. In the experiments, we first load the data into memory and then build the sketches. The computation time reported does not include computing the final approximation $[[\hat{\bm{A}}]]_r$ since in practice it is usually stored in the factorized form for subsequent tasks. The need to compute $[[\hat{\bm{A}}]]_r$ is only for the error metric. We also provide visual comparisons on the first few left singular vectors in each case.

\subsection{Yale Face Dataset} 

The data is originally from the Yale Face Database B \cite{Georghiades2001}. For this dataset, $\bm{A}\in\mathbb{R}^{2500\times 640}$ since there are 640 face images of size $50\times 50$. Based on the scree plot, we choose $r=20$. The optimal $err$ is 0.033. 

\begin{figure}[H]
\label{fig:scree_Yale}
    \centering
    \includegraphics[width=0.35\textwidth]{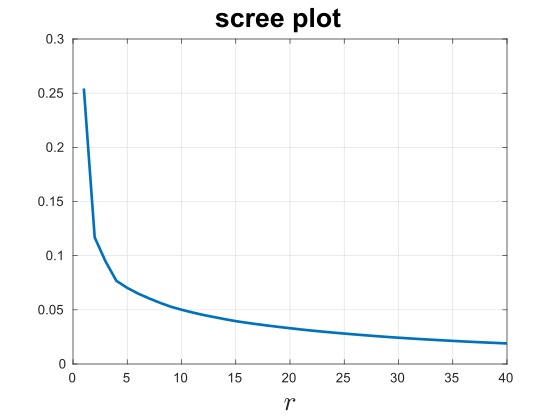}
    \caption{Scree plot for Yale Face dataset.}
\end{figure}

\begin{table}[H]
\begin{center}
\small
\begin{tabular}{ccccc}
\hline
\multicolumn{1}{|c|}{} & \multicolumn{1}{c|}{SketchySVD} & \multicolumn{3}{c|}{SketchyCoreSVD} \\
\hline
\multicolumn{1}{|c|}{$p$} & \multicolumn{1}{c|}{-} & \multicolumn{1}{c|}{0.3} & \multicolumn{1}{c|}{0.35} & \multicolumn{1}{c|}{0.4} \\
\hline
\multicolumn{1}{|c|}{$err$} & \multicolumn{1}{c|}{0.066} & \multicolumn{1}{c|}{0.0765} & \multicolumn{1}{c|}{0.0737} & \multicolumn{1}{c|}{0.0717} \\
\hline
\multicolumn{1}{|c|}{time (sec)} & \multicolumn{1}{c|}{0.0239} & \multicolumn{1}{c|}{0.0134} & \multicolumn{1}{c|}{0.0146} & \multicolumn{1}{c|}{0.0166}  \\
\hline
\end{tabular}
\caption{Performance comparisons for Yale Face dataset.}
\end{center}
\end{table}

We choose $k=4r+1=81$ and $s=2k+1=163$ for both methods. In Table $1$, we see that the error of SketchySVD is about twice the optimal error. Our SketchyCoreSVD can achieve about the same error bound with less time . As the sampling ratio $p$ increases from 30\% to 40\%, the $err$ decreases while computation time only gradually increases. In the visual comparison Figure $2$, we show that the singular vector(s) can be estimated accurately, and with sampling ratio $p=40\%$ for SketchyCoreSVD.

\begin{figure}[H]
    \centering
    \includegraphics[width=0.3\textwidth]{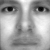}
    \includegraphics[width=0.3\textwidth]{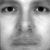}
    \includegraphics[width=0.3\textwidth]{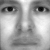}
    \caption{Visual comparison for the first left singular vector computed on the Yale Face dataset. Ground truth (left column), SketchySVD (middle column), SketchyCoreSVD (right column, $p=0.4$). The peak signal to noise ratio (PSNR) is 44.0320 for SketchySVD, and 44.2639 for SketchyCoreSVD.}
\end{figure}

\subsection{Cardiac Magnetic Resonance Imaging}

We have a collection of time-varying, 2D slice stacks of  cardiac magnetic resonance images (Cardiac MRI). We select a 2D MR spatial snapshot and consider a time-sequence of this in the form of a  data matrix  $\bm{A}\in\mathbb{R}^{45056\times 160}$  (i.e., 160 time snapshots 2D images, each of of size $256\times 176$). Based on the scree plot, we choose $r=5$. The optimal $err$ is 0.0011.

\begin{figure}[H]
\label{fig:scree_Cardiac}
    \centering
    \includegraphics[width=0.35\textwidth]{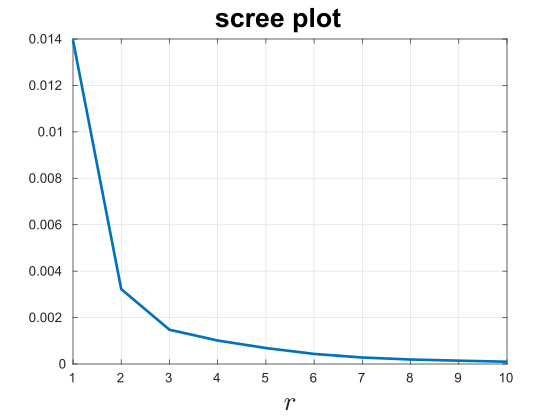}
    \caption{Scree plot for Cardiac MRI dataset.}
\end{figure}

We choose $k=4r+1=21$ and $s=2k+1=43$ for both methods. In Table $2$, we see that the error of SketchySVD is less than twice the optimal error. Our SketchyCoreSVD can achieve the same error bound in less time. As sampling ratio $p$ increases from 30\% to 40\%, $err$ decreases while computation time gradually increases. In the visual comparison Figure $4$, we show that the singular vector(s) can be estimated accurately, and with sampling ratio $p=40\%$ for SketchyCoreSVD.

\begin{table}[H]
\begin{center}
\small
\begin{tabular}{ccccc}
\hline
\multicolumn{1}{|c|}{} & \multicolumn{1}{c|}{SketchySVD} & \multicolumn{3}{c|}{SketchyCoreSVD} \\
\hline
\multicolumn{1}{|c|}{$p$} & \multicolumn{1}{c|}{-} & \multicolumn{1}{c|}{0.3} & \multicolumn{1}{c|}{0.35} & \multicolumn{1}{c|}{0.4} \\
\hline
\multicolumn{1}{|c|}{$err$} & \multicolumn{1}{c|}{0.0019} & \multicolumn{1}{c|}{0.0021} & \multicolumn{1}{c|}{0.0021} & \multicolumn{1}{c|}{0.0019} \\
\hline
\multicolumn{1}{|c|}{time (sec)} & \multicolumn{1}{c|}{0.0567} & \multicolumn{1}{c|}{0.0316} & \multicolumn{1}{c|}{0.038} & \multicolumn{1}{c|}{0.0396}  \\
\hline
\end{tabular}
\caption{Performance comparisons for Cardiac MRI dataset.}
\end{center}
\end{table}

\begin{figure}[H]
    \centering
    \includegraphics[width=0.3\textwidth]{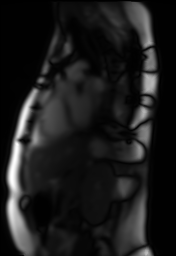}
    \includegraphics[width=0.3\textwidth]{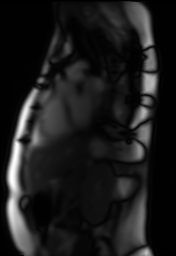}
    \includegraphics[width=0.3\textwidth]{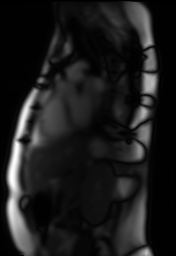}\\
    \includegraphics[width=0.3\textwidth]{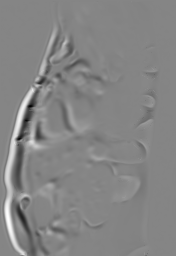}
    \includegraphics[width=0.3\textwidth]{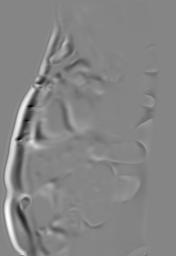}
    \includegraphics[width=0.3\textwidth]{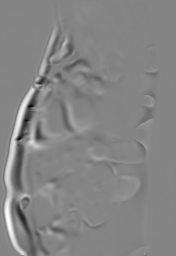}\\
    \caption{Visual comparison for the first two left singular vectors computed on the Cardiac MRI dataset. Ground truth (left column), SketchySVD (middle column), SketchyCoreSVD (right column, $p=0.4$). The PSNR values are $\{68.5232,59.1094\}$ for SketchySVD, and $\{69.5770,58.3338\}$ for SketchyCoreSVD.}
\end{figure}

\subsection{The BR1003 Breast Cancer Dataset}

This breast cancer data matrix we use is  extracted from a collection of Fourier-transform infrared spectroscopy (FTIR) spectral signatures of breast tissues \cite{Plos-bajaj2019}.  The matrix $\bm{A}\in\mathbb{R}^{1506\times 783090}$ is created from 789030 spectral signatures of length 1506. Based on the scree plot, we choose $r=6$. The optimal $err$ is 0.002. 

\begin{figure}[H]
\label{fig:scree_BR1003}
    \centering
    \includegraphics[width=0.35\textwidth]{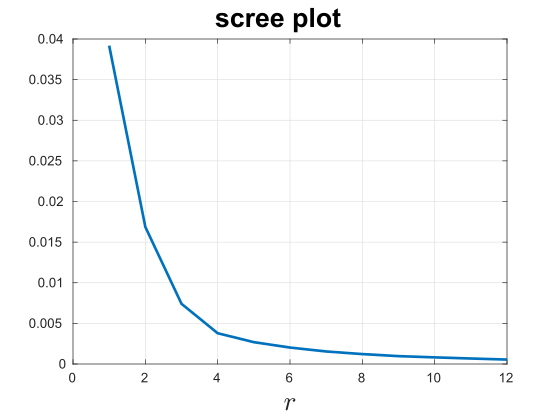}
    \caption{Scree plot for BR1003 dataset.}
\end{figure}

\begin{table}[H]
\begin{center}
\small
\begin{tabular}{ccccc}
\hline
\multicolumn{1}{|c|}{} & \multicolumn{1}{c|}{SketchySVD} & \multicolumn{3}{c|}{SketchyCoreSVD} \\
\hline
\multicolumn{1}{|c|}{$p$} & \multicolumn{1}{c|}{-} & \multicolumn{1}{c|}{0.04} & \multicolumn{1}{c|}{0.06} & \multicolumn{1}{c|}{0.08} \\
\hline
\multicolumn{1}{|c|}{$err$} & \multicolumn{1}{c|}{0.0025} & \multicolumn{1}{c|}{0.0031} & \multicolumn{1}{c|}{0.0029} & \multicolumn{1}{c|}{0.0027} \\
\hline
\multicolumn{1}{|c|}{time (sec)} & \multicolumn{1}{c|}{3.2049} & \multicolumn{1}{c|}{0.816} & \multicolumn{1}{c|}{1.1643} & \multicolumn{1}{c|}{1.3614}  \\
\hline
\end{tabular}
\caption{Performance comparisons for BR1003 dataset.}
\end{center}
\end{table}

\begin{figure}[H]
    \centering
    \includegraphics[width=0.3\textwidth]{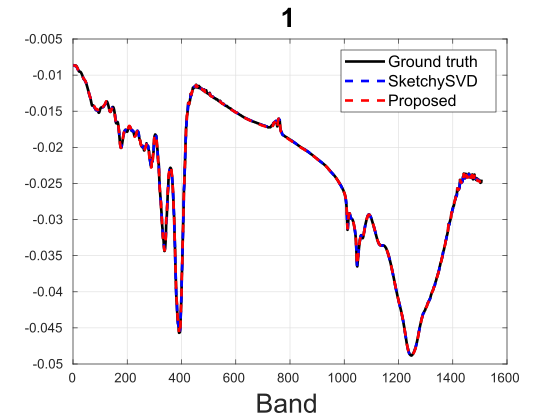}
    \includegraphics[width=0.3\textwidth]{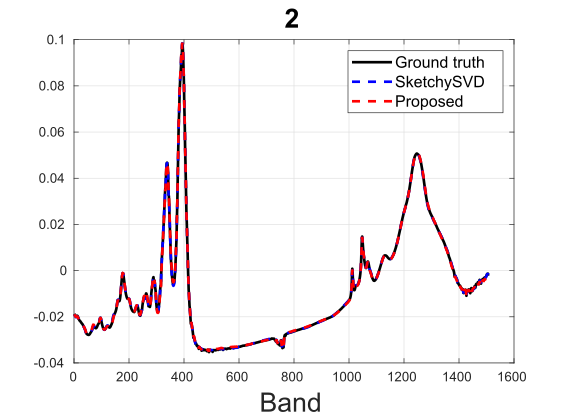}
    \includegraphics[width=0.3\textwidth]{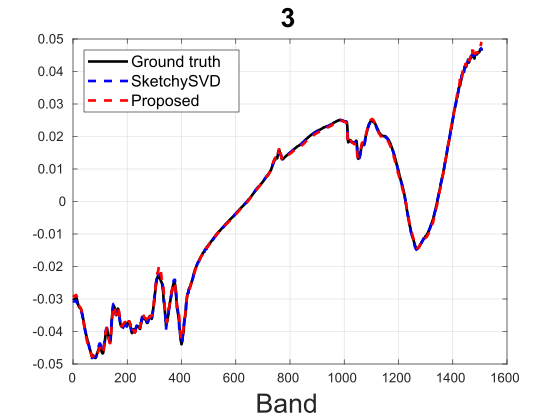}\\
    \includegraphics[width=0.3\textwidth]{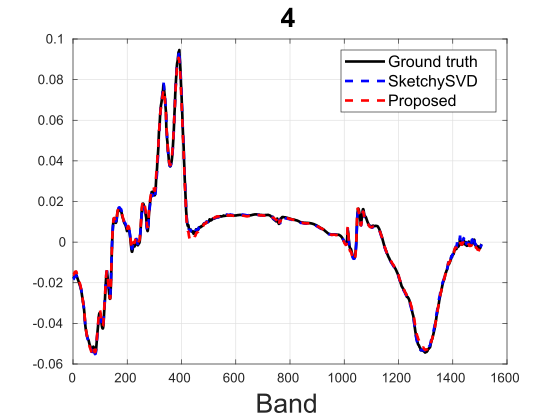}
    \includegraphics[width=0.3\textwidth]{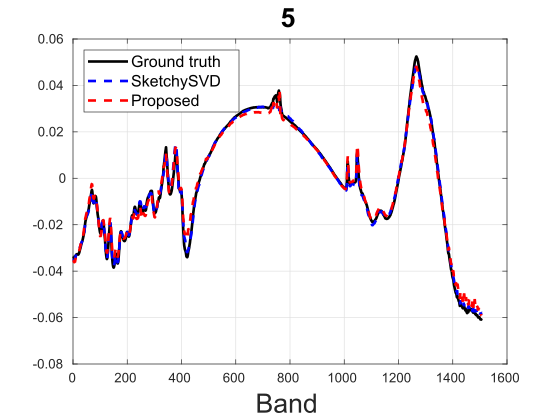}
    \includegraphics[width=0.3\textwidth]{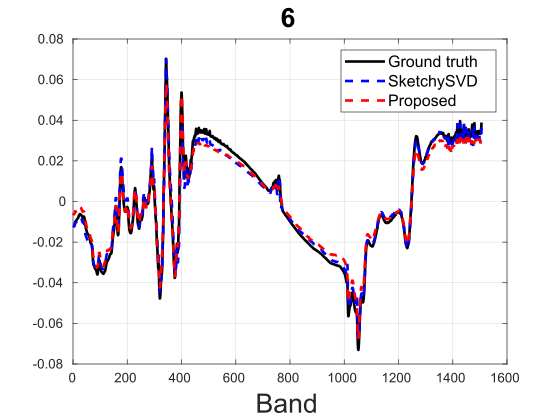}
    \caption{Visual comparison for the first six left singular vectors computed on the BR1003 dataset, where $p=0.08$ for our SketchyCoreSVD method.}
\end{figure}

\begin{table}[H]
\begin{center}
\small
\begin{tabular}{ccccccc}
\hline
\multicolumn{1}{|c|}{} & \multicolumn{1}{c|}{1} & \multicolumn{1}{c|}{2} & \multicolumn{1}{c|}{3} & \multicolumn{1}{c|}{4} & \multicolumn{1}{c|}{5} & \multicolumn{1}{c|}{6} \\
\hline
\multicolumn{1}{|c|}{SketchySVD} & \multicolumn{1}{c|}{56.6010} & \multicolumn{1}{c|}{41.9409} & \multicolumn{1}{c|}{37.9187} & \multicolumn{1}{c|}{32.6139} & \multicolumn{1}{c|}{25.3890} & \multicolumn{1}{c|}{18.6244} \\
\hline
\multicolumn{1}{|c|}{SketchyCoreSVD} & \multicolumn{1}{c|}{54.2577} & \multicolumn{1}{c|}{39.0797} & \multicolumn{1}{c|}{34.9353} & \multicolumn{1}{c|}{31.8365} & \multicolumn{1}{c|}{23.0102} & \multicolumn{1}{c|}{15.6122} \\
\hline
\end{tabular}
\caption{Comparison of signal to noise (SNR) ratios for the first six left singular vectors computed on the BR1003 dataset, where $p=0.08$ for our SketchyCoreSVD method.}
\end{center}
\end{table}

We choose $k=4r+1=25$ and $s=2k+1=51$ for both methods. In Table $3$, we see that the error of SketchySVD is about the same as the optimal error. SketchyCoreSVD can achieve about the same error bound in less than half of the time. As the sampling ratio $p$ increases from 4\% to 8\%, the $err$ decreases while the computation time only gradually increases. In the visual comparison Figure $6$, we show that the singular vector(s) can be estimated accurately, and with sampling ratio $p=8\%$ for SketchyCoreSVD.

\subsection{Video Dataset}

This dataset is a color video of size $1080\times 1920\times 3\times 2498$. It was originally used by \cite{malik2018low} to test tensor approximations. We converted the video into to grayscale, reduced the spatial size by a factor of 2, and discarded the first 100 and last 198 frames due to camera blur. The resulting data matrix $\bm{A}\in\mathbb{R}^{518400\times 2200}$. Based on the scree plot, we choose $r=25$. The optimal $err$ is 0.0066.

\begin{figure}[H]
\label{fig:scree_Video}
    \centering
    \includegraphics[width=0.35\textwidth]{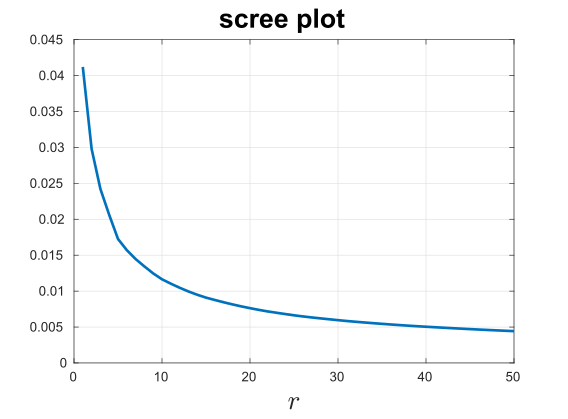}
    \caption{Scree plot for Video dataset.}
\end{figure}

\begin{table}[H]
\begin{center}
\small
\begin{tabular}{ccccc}
\hline
\multicolumn{1}{|c|}{} & \multicolumn{1}{c|}{SketchySVD} & \multicolumn{3}{c|}{SketchyCoreSVD} \\
\hline
\multicolumn{1}{|c|}{$p$} & \multicolumn{1}{c|}{-} & \multicolumn{1}{c|}{0.1} & \multicolumn{1}{c|}{0.15} & \multicolumn{1}{c|}{0.2} \\
\hline
\multicolumn{1}{|c|}{$err$} & \multicolumn{1}{c|}{0.0148} & \multicolumn{1}{c|}{0.0213} & \multicolumn{1}{c|}{0.0177} & \multicolumn{1}{c|}{0.0165} \\
\hline
\multicolumn{1}{|c|}{time (sec)} & \multicolumn{1}{c|}{8.0062} & \multicolumn{1}{c|}{2.3224} & \multicolumn{1}{c|}{2.6266} & \multicolumn{1}{c|}{3.5619}  \\
\hline
\end{tabular}
\caption{Performance comparisons for Video dataset.}
\end{center}
\end{table}

We choose $k=4r+1=101$ and $s=2k+1=203$ for both methods. In Table $4$, we see that the error of SketchySVD is about twice the optimal error. SketchyCoreSVD can achieve about the same error bound in less than half of the time. As the sampling ratio $p$ is  increased from 10\% to 20\%, the $err$ decreases while the computation time only gradually increases. In the visual comparison Figure $8$, we show that the singular vector(s) can be estimated accurately, and with sampling ratio $p=20\%$ for SketchyCoreSVD.

\begin{figure}[H]
    \centering
    \includegraphics[width=0.3\textwidth]{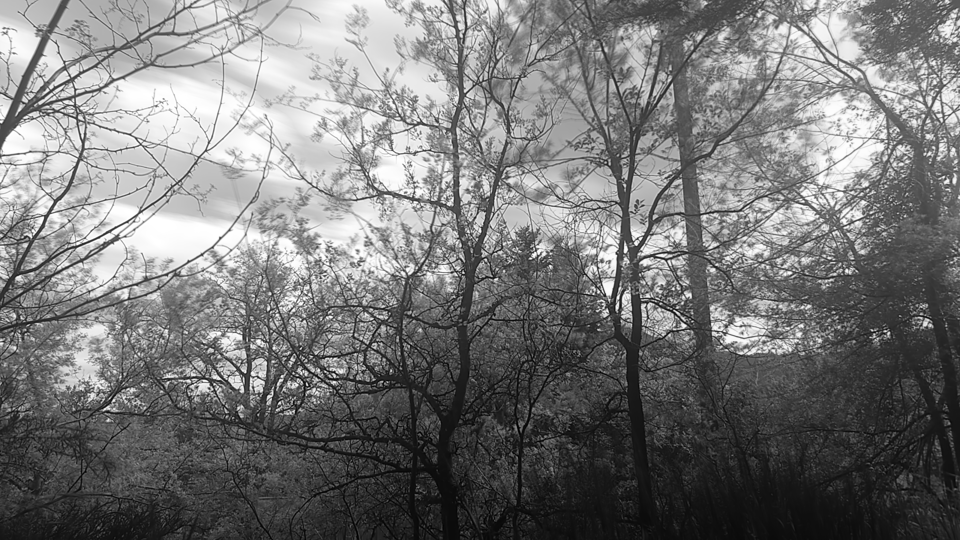}
    \includegraphics[width=0.3\textwidth]{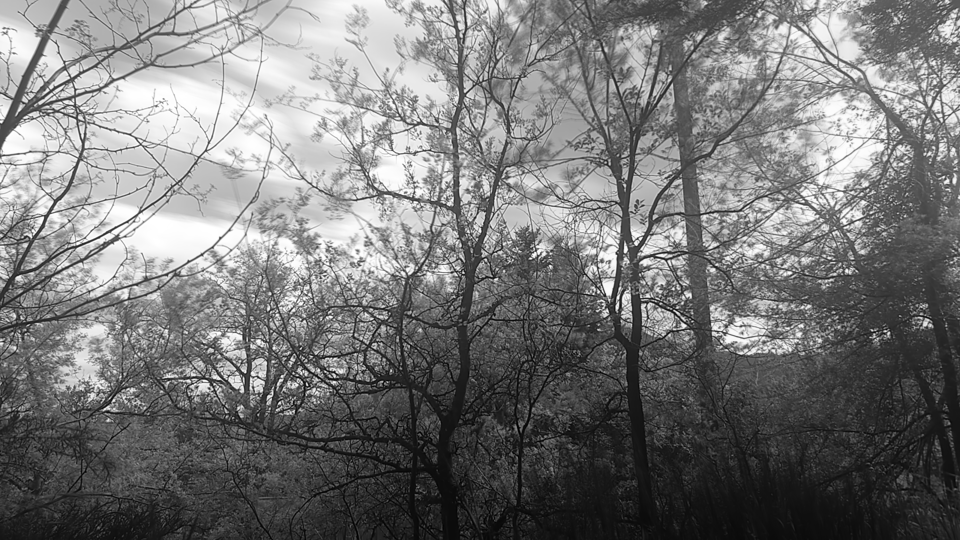}
    \includegraphics[width=0.3\textwidth]{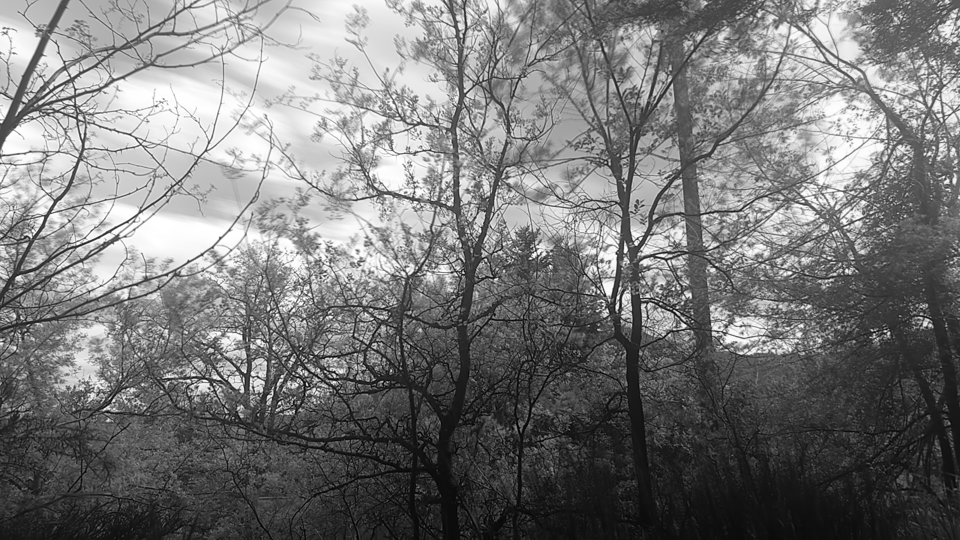}\\
    \includegraphics[width=0.3\textwidth]{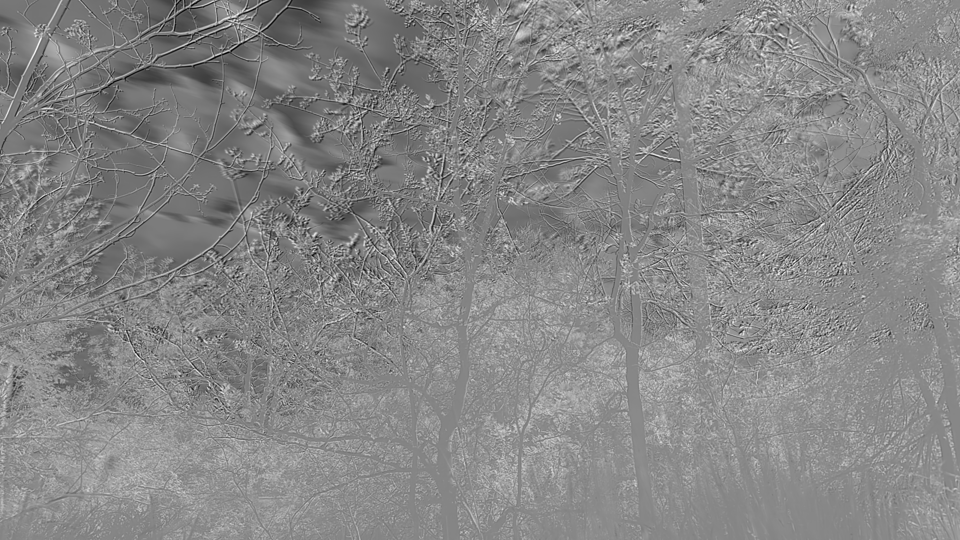}
    \includegraphics[width=0.3\textwidth]{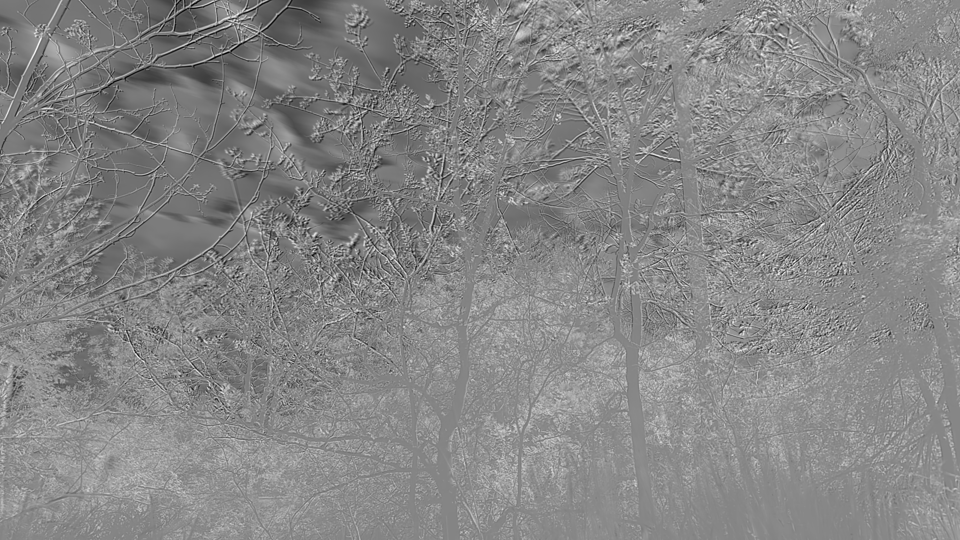}
    \includegraphics[width=0.3\textwidth]{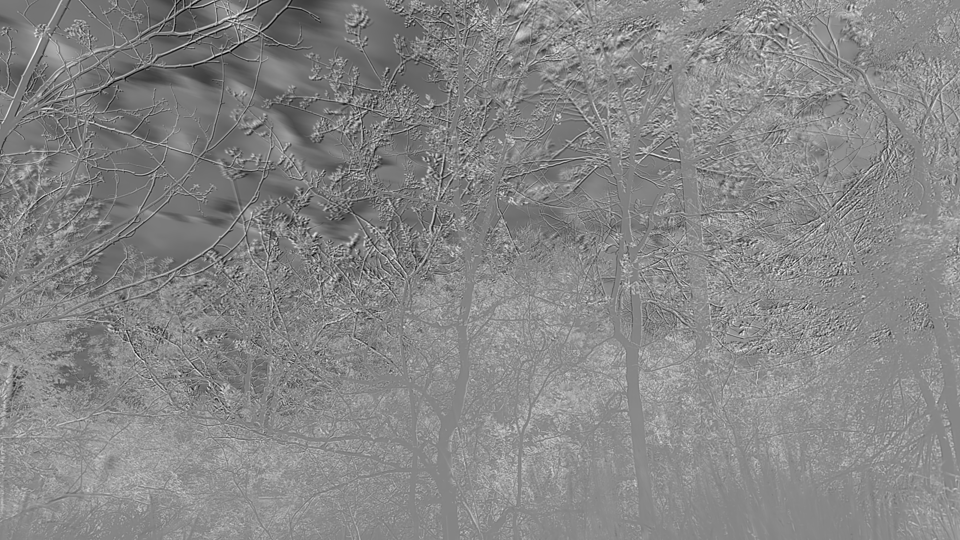}\\
    \caption{Visual comparison for the first two left singular vectors computed on the Video dataset. Ground truth (left column), SketchySVD (middle column), SketchyCoreSVD (right column, $p=0.2$). The PSNR values are $\{55.5701,50.2379\}$ for SketchySVD, and $\{55.8322,49.2919\}$ for SketchyCoreSVD.}
\end{figure} 

\noindent \textbf{Acknowledgment} 
This research was supported in part by a grant from NIH - R01GM117594.

\bibliographystyle{plain}
\bibliography{references}

\appendix
\section{Algorithms Based on Sketches}
\label{sec:prior_algorithms}

We provide a brief comparison of prior algorithms based on sketches. 

\fbox{\begin{minipage}{15cm}
\noindent\textbf{\cite{halko2011finding}} builds sketches
$$
\bm{X} = \bm{\Gamma}\bm{A}
\quad\text{and}\quad
\bm{Y} = \bm{A}\bm{\Omega}^*.
$$
Compute $(\bm{P},\sim,\sim)=\text{svds}(\bm{X}^*,r)$, $(\bm{Q},\sim,\sim)=\text{svds}(\bm{Y},r)$; 
$\bm{C}_1 = \bm{Q}^*\bm{Y}\left((\bm{\Omega}\bm{P})^{\dagger}\right)^*$, $\bm{C}_2 = (\bm{\Gamma}\bm{Q})^{\dagger}\bm{X}\bm{P}$;
$(\bm{U},\bm{\Sigma},\bm{V})=\text{svd}((\bm{C}_1+\bm{C}_2)/2)$;
$$
\hat{\bm{A}} = (\bm{Q}\bm{U})\bm{\Sigma}(\bm{P}\bm{V})^*.
$$
\end{minipage}}

\fbox{\begin{minipage}{15cm}
\noindent\textbf{\cite{woodruff2014sketching}} builds sketches
$$
\bm{X} = \bm{\Gamma}\bm{A}
\quad\text{and}\quad
\bm{Y} = \bm{A}\bm{\Omega}^*.
$$
Compute $\bm{Y}=\bm{Q}\bm{R}$; $\bm{\Gamma}\bm{Q}=\bm{U}\bm{T}$;
$$
\hat{\bm{A}} = \bm{Q}\bm{T}^{\dagger}[[\bm{U}^*\bm{Y}]]_r.
$$
\end{minipage}}

\fbox{\begin{minipage}{15cm}
\noindent\textbf{\cite{cohen2015dimensionality}} builds sketches
$$
\bm{X} = \bm{\Gamma}\bm{A}
\quad\text{and}\quad
\bm{Y} = \bm{A}\bm{\Omega}^*.
$$
Compute $(\bm{V},\sim,\sim)=\text{svds}(\bm{Y},r)$; $\bm{\Gamma}\bm{Q}=\bm{U}\bm{T}$;
$$
\hat{\bm{A}} = \bm{V}\bm{T}^{\dagger}[[\bm{U}^*\bm{Y}]]_r.
$$
\end{minipage}}

\fbox{\begin{minipage}{15cm}
\noindent\textbf{\cite{boutsidis2016optimal}} with simplification suggested by \cite{upadhyay2018price}. Build sketches
$$
\bm{X} = \bm{\Gamma}\bm{A},
\quad
\bm{Y} = \bm{A}\bm{\Omega}^*,
\quad\text{and}\quad
\bm{Z} = \bm{\Phi}\bm{A}\bm{\Psi}^*.
$$
Compute $\bm{X}^*=\bm{P}\bm{R}_1$, $\bm{Y}=\bm{Q}\bm{R}_2$; $\bm{\Phi}\bm{Q}=\bm{U}_1\bm{\Sigma}_1\bm{V}_1^*$, $\bm{\Psi}\bm{P}=\bm{U}_2\bm{\Sigma}_2\bm{V}_2^*$;
$$
\hat{\bm{A}} = \bm{Q}\bm{V}_1\bm{\Sigma}_1^{\dagger}[[\bm{U}_1^*\bm{Z}\bm{U}_2]]_r\bm{\Sigma}_2^{\dagger}\bm{V}_2^*\bm{P}^*.
$$
\end{minipage}}

\fbox{\begin{minipage}{15cm}
\noindent\textbf{\cite{tropp2017practical}} builds sketches
$$
\bm{X} = \bm{\Gamma}\bm{A}
\quad\text{and}\quad
\bm{Y} = \bm{A}\bm{\Omega}^*.
$$
Compute $\bm{Y}=\bm{Q}\bm{R}$;
$$
\hat{\bm{A}} = \bm{Q}[[(\bm{\Gamma}\bm{Q})^{\dagger}\bm{X}]]_r.
$$
\end{minipage}} 
\section{Proof of Theorem 1}

The following proofs are for the approximation of the column space of $\bm{A}$, i.e., $\|\bm{A}-\bm{Q}\bm{Q}^*\bm{A}\|_F$. Similar argument holds for $\|\bm{A}-\bm{A}\bm{P}\bm{P}^*\|_F$. 

We first state a few auxiliary lemmas to facilitate the presentation of the main proof of Theorem 1. Lemma~\ref{lem:incoherence}, adapted from \cite[Lemma~3.4]{tropp2011improved}, is essential for our proofs since it provides a lower bound on the singular values of a submatrix from randomly sampled rows of an ``incoherent'' orthonormal matrix. 

\begin{lemma}\label{lem:incoherence}
Suppose $\max_{j}\{\|\bm{V}_1^{(j,:)}\|_2\}\leq \sqrt{\frac{\mu r}{N}}$. For $\alpha>0$, select the sample size
$$
n\geq 8\mu r\log r. 
$$
Draw a random subset $\Theta$ from $\{1,\cdots,N\}$ by sampling $n$ coordinates without replacement. Then with probability at least $1-\frac{2}{r^3}$,
$$
\sqrt{\frac{n}{6N}}\leq \sigma_{r}(\bm{V}_1^{(\Theta,:)})\quad
\text{and}\quad
\sigma_1(\bm{V}_1^{(\Theta,:)}) \leq \sqrt{\frac{13n}{6N}}.
$$ 
\end{lemma}

If we further apply Gaussian maps to the subsampled rows, Lemma \ref{lem:full_rank} ensures that the full rankness is preserved almost surely. 
\begin{lemma}\label{lem:full_rank}
Suppose $\bm{V}_1^{(\Theta,:)}\in\mathbb{R}^{n\times r}$ is of full column rank, and $\bm{\Omega}\in\mathbb{R}^{k\times n}$ is standard normal Gaussian. Then
$\bm{\Omega}\bm{V}_1^{(\Theta,:)}\in\mathbb{R}^{k\times r}$ is of full column rank almost surely.
\end{lemma}
\begin{proof}
We prove for $\bm{\Omega}\bm{V}_1^{(\Theta,:)}$. From Lemma~\ref{lem:incoherence}, we know that $\bm{V}_1^{(\Theta,:)}$ is of full column rank. Denote its truncated SVD by $\bm{W}_1\bm{\Lambda}\bm{W}_2^*$, where $\bm{W}_1\in\mathbb{R}^{n\times r}$, $\bm{\Lambda}\in\mathbb{R}^{r\times r}$, and $\bm{W}_2\in\mathbb{R}^{r\times r}$. Since $\bm{W}_1$ is orthonormal, $\bm{\Omega}\bm{W}_1\in\mathbb{R}^{k\times r}$ is standard Gaussian. Thus with probability one, $\bm{\Omega}\bm{W}_1$ is of full column rank. Consequently, $\bm{\Omega}\bm{V}_1^{(\Theta,:)}=(\bm{\Omega}\bm{W}_1)\bm{\Lambda}\bm{W}_2^*$ is also of full column 
rank.
\end{proof}

From Lemma~\ref{lem:full_rank}, we know that $(\bm{V}_1^*)^{(:,\Theta)}\bm{\Omega}^*\in\mathbb{R}^{r\times k}$ is of full row rank. Define $\bm{\Omega}_1=(\bm{V}_1^*)^{(:,\Theta)}\bm{\Omega}^*$, and $\bm{\Omega}_2=(\bm{V}_2^*)^{(:,\Theta)}\bm{\Omega}^*$. Lemma~\ref{lem:deterministic} below provides a deterministic bound of $\|\bm{A}-\bm{Q}\bm{Q}^*\bm{A}\|_F^2$. It can be proved similarly as \cite[Theorem~9.1]{halko2011finding}. For completeness, we include its proof in Section~\ref{sec:Lemma3}.
\begin{lemma}\label{lem:deterministic}
Assuming that $\bm{\Omega}_1\in\mathbb{R}^{r\times k}$ has full row rank, and $\bm{Q}\in\mathbb{R}^{M\times k}$ is the orthonormal columns from the QR decomposition of $\bm{A}^{(:,\Theta)}\bm{\Omega}$. Then
\begin{align}\label{eq:deterministic}
\|\bm{A}-\bm{Q}\bm{Q}^*\bm{A}\|_F^2 \leq \|\bm{\Sigma}_2\|_F^2+\|\bm{\Sigma}_2\bm{\Omega}_2\bm{\Omega}_1^{\dagger}\|_F^2,
\end{align}
where $\bm{\Sigma}_2$ is the diagonal matrix containing the $(r+1)$ to $\min\{M,N\}$ singular values of $\bm{A}$.
\end{lemma}

For the bound in \eqref{eq:deterministic}, $\bm{\Omega}_2$ and $\bm{\Omega}_1$ are not independent in our case, contrary to what \cite{tropp2019streaming} deals with. Our strategy is to separately bound $\bm{\Sigma}_2\bm{\Omega}_2$ and $\bm{\Omega}_1^{\dagger}$. For $\bm{\Omega}_1^{\dagger}$, we need to deal with the operator norm of the pseudo-inverse of random Gaussian at some point. Lemma~\ref{lem:pseudo_inverse}, adapted from \cite[Proposition~10.3]{halko2011finding}, provides such a bound. 

\begin{lemma}\label{lem:pseudo_inverse}
For $\bm{G}\in\mathbb{R}^{k\times r}$ ($k\geq r+4$) being standard normal Gaussian,
$$
\mathbb{P}\left\{\|\bm{G}^{\dagger}\|_2\geq \frac{e\sqrt{k}}{k-r+1}\cdot t\right\}\leq t^{-(k-r+1)},\quad \forall t\geq 1.
$$
Taking $t=k^{\frac{3}{k-r+1}}$, then 
$$
\|\bm{G}^{\dagger}\|_2\leq\frac{e\sqrt{k}}{k-r+1}\cdot k^{\frac{3}{k-r+1}}
$$ 
with probability at least $1-\frac{1}{k^3}$.
\end{lemma}

We bound $\|\bm{\Sigma}_2\bm{\Omega}_2\|_F$ with the help of Lemma~\ref{lem:expectation} and Lemma~\ref{lem:tail}.

\begin{lemma}{(\cite[Proposition~10.1]{halko2011finding})}\label{lem:expectation}
Fix matrices $\bm{S}$, $\bm{T}$, and draw a standard Gaussian matrix $\bm{G}$. Then
$$
\mathbb{E}(\|\bm{S}\bm{G}\bm{T}\|_F^2) = \|\bm{S}\|_F^2\|\bm{T}\|_F^2.
$$
\end{lemma}

\begin{lemma}{(\cite[Proposition~10.3]{halko2011finding})}\label{lem:tail}
Suppose $h$ is a Lipschitz function on matrices:
$$
|h(\bm{X})-h(\bm{Y})|\leq L\|\bm{X}-\bm{Y}\|_F,\quad\forall\bm{X},\bm{Y}.
$$
Draw a standard Gaussian matrix $\bm{G}$. Then
$$
\mathbb{P}\{h(\bm{G})\geq \mathbb{E}(h(\bm{G}))+Lt\}\leq e^{-t^2/2}.
$$
\end{lemma}

Now we are ready to present the main proof of Theorem 1.

\begin{align*}
\|\bm{A}-\bm{Q}\bm{Q}^*\bm{A}\|_F^2 
\stackrel{(a)}{\leq}&\|\bm{\Sigma}_2\|_F^2+\|\bm{\Sigma}_2\bm{\Omega}_2\bm{\Omega}_1^{\dagger}\|_F^2\\
= &\|\bm{\Sigma}_2\|_F^2+\|\bm{\Sigma}_2(\bm{V}_2^*)^{(:,\Theta)}\bm{\Omega}^*\bm{\Omega}_1^{\dagger}\|_F^2\\
= &\|\bm{\Sigma}_2\|_F^2+\|\bm{\Sigma}_2(\bm{V}_2^*)^{(:,\Theta)}\bm{\Omega}^*((\bm{V}_1^*)^{(:,\Theta)}\bm{\Omega}^*)^{\dagger}\|_F^2\\
\stackrel{(b)}{=} &\|\bm{\Sigma}_2\|_F^2+\|\bm{\Sigma}_2(\bm{V}_2^*)^{(:,\Theta)}\bm{\Omega}^*(\bm{W}_2\bm{\Lambda}\bm{W}_1^*\bm{\Omega}^*)^{\dagger}\|_F^2\\
\stackrel{(c)}{=} &\|\bm{\Sigma}_2\|_F^2+\|\bm{\Sigma}_2(\bm{V}_2^*)^{(:,\Theta)}\bm{\Omega}^*(\bm{W}_1^*\bm{\Omega}^*)^{\dagger}\bm{\Lambda}^{\dagger}\bm{W}_2^{\dagger}\|_F^2\\
\leq &\|\bm{\Sigma}_2\|_F^2+\|(\bm{W}_1^*\bm{\Omega}^*)^{\dagger}\bm{\Lambda}^{\dagger}\bm{W}_2^{\dagger}\|_2^2\cdot\|\bm{\Sigma}_2(\bm{V}_2^*)^{(:,\Theta)}\bm{\Omega}^*\|_F^2\\
\leq &\|\bm{\Sigma}_2\|_F^2+\|\bm{\Lambda}^{\dagger}\|_2^2\cdot\|(\bm{W}_1^*\bm{\Omega}^*)^{\dagger}\|_2^2\cdot\|\bm{\Sigma}_2(\bm{V}_2^*)^{(:,\Theta)}\bm{\Omega}^*\|_F^2\\
\stackrel{(d)}{\leq} &\|\bm{\Sigma}_2\|_F^2+\frac{6N}{n}\cdot\|(\bm{W}_1^*\bm{\Omega}^*)^{\dagger}\|_2^2\cdot\|\bm{\Sigma}_2(\bm{V}_2^*)^{(:,\Theta)}\bm{\Omega}^*\|_F^2\\
\stackrel{(e)}{\leq} &\|\bm{\Sigma}_2\|_F^2+\frac{6N}{n}\cdot\frac{e^2k^{1+6/(k-r+1)}}{(k-r+1)^2} \cdot\underbrace{\|\bm{\Sigma}_2(\bm{V}_2^*)^{(:,\Theta)}\bm{\Omega}^*\|_F^2}_{T_1},
\end{align*}
where (a) holds with probability at least $1-\frac{2}{r^3}$, by applying Lemma~\ref{lem:incoherence}, Lemma~\ref{lem:full_rank}, and Lemma~\ref{lem:deterministic}; in $(b)$, we assume $\bm{V}_1^{(\Theta,:)}=\bm{W}_1\bm{\Lambda}\bm{W}_2^*$ is the SVD (of rank $r$); $(c)$ is due to the property of pseudo-inverse; $(d)$ is due to the lower bound of the singular value in Lemma 1; (e) holds with probability at least $1-\frac{1}{k^3}$, by applying Lemma 4 to $\bm{\Omega}\bm{W}_1$ .

\textbf{Bound for $T_1$.} By Lemma~\ref{lem:expectation},
$$
\mathbb{E}(\|\bm{\Sigma}_2(\bm{V}_2^*)^{(:,\Theta)}\bm{\Omega}^*\|_F^2) = \|\bm{\Sigma}_2(\bm{V}_2^*)^{(:,\Theta)}\|_F^2\cdot\|\bm{I}_k\|_F^2\leq k\|\bm{\Sigma}_2\|_F^2,
$$
where $\bm{I}_k$ is the identity matrix of size $k$, and the last inequality is due to $\|(\bm{V}_2^*)^{(:,\Theta)}\|_2=\|\bm{V}_2^{(\Theta,:)}\|_2\leq 1$. Consider the function $h(\bm{X})=\|\bm{\Sigma}_2(\bm{V}_2^*)^{(:,\Theta)}\bm{X}\|_F$,
\begin{align*}
|h(\bm{X})-h(\bm{Y})|=&|\|\bm{\Sigma}_2(\bm{V}_2^*)^{(:,\Theta)}\bm{X}\|_F-\|\bm{\Sigma}_2(\bm{V}_2^*)^{(:,\Theta)}\bm{Y}\|_F| \\
\leq&\|\bm{\Sigma}_2(\bm{V}_2^*)^{(:,\Theta)}(\bm{X}-\bm{Y})\|_F\\
\leq&\|\bm{\Sigma}_2(\bm{V}_2^*)^{(:,\Theta)}\|_2\cdot\|\bm{X}-\bm{Y}\|_F\leq\|\bm{\Sigma}_2\|_2\cdot\|\bm{X}-\bm{Y}\|_F.
\end{align*}
Thus the Lipschitz constant is $\|\bm{\Sigma}_2\|_2$. Now by Lemma~\ref{lem:tail}, with probability at least $1-e^{-t^2/2}$,
$$
\|\bm{\Sigma}_2(\bm{V}_2^*)^{(:,\Theta)}\bm{\Omega}^*\|_F\leq \sqrt{k}\|\bm{\Sigma}_2\|_F+\|\bm{\Sigma}_2\|_2\cdot t.
$$
A particular choice of $t=\sqrt{6\log k}$ gives 
$$
\|\bm{\Sigma}_2(\bm{V}_2^*)^{(:,\Theta)}\bm{\Omega}^*\|_F\leq \sqrt{k}\|\bm{\Sigma}_2\|_F+\sqrt{6\log k}\|\bm{\Sigma}_2\|_2
$$ 
with probability at least $1-\frac{1}{k^3}$. 

Based on the derived bound for $T_1$,
\begin{align*}
\|\bm{A}-\bm{Q}\bm{Q}^*\bm{A}\|_F
\leq&\sqrt{\|\bm{\Sigma}_2\|_F^2+\frac{6N}{n}\cdot\frac{e^2k^{1+6/(k-r+1)}}{(k-r+1)^2} \cdot\left(\sqrt{k}\|\bm{\Sigma}_2\|_F+\sqrt{6\log k}\|\bm{\Sigma}_2\|_2\right)^2}\\
\leq &\|\bm{\Sigma}_2\|_F+\sqrt{\frac{6N}{n}}\cdot\frac{e\sqrt{k}}{k-r+1}\cdot k^{\frac{3}{k-r+1}}\cdot\left(\sqrt{k}\|\bm{\Sigma}_2\|_F+\sqrt{6\log k}\|\bm{\Sigma}_2\|_2\right)\\
= & (C_1(p,k,r)+1)\cdot\|\bm{\Sigma}_2\|_F+C_2(p,k,r)\cdot\|\bm{\Sigma}_2\|_2,
\end{align*}
where the second inequality is due to $\sqrt{a^2+b^2}\leq a+b$ for $a,b\geq 0$, and in the last equality,
$$
C_1(p,k,r)=\sqrt{\frac{6e^2}{p}}\cdot\frac{k}{k-r+1}\cdot k^{\frac{3}{k-r+1}},
$$
$$
C_2(p,k,r)=\sqrt{\frac{36e^2}{p}}\cdot\frac{\sqrt{k\log k}}{k-r+1}\cdot k^{\frac{3}{k-r+1}}.
$$

\subsection{Proof of Lemma~\ref{lem:deterministic}}\label{sec:Lemma3}
First, the left unitary factor $\bm{U}$ plays no essential role. To see this, define
$$
\widetilde{\bm{A}} =  \bm{U}^*\bm{A} =
\left[\begin{array}{c}
\bm{\Sigma}_1\bm{V}_1^* \\
\bm{\Sigma}_2\bm{V}_2^*
\end{array}\right]\quad
\text{and}\quad
\widetilde{\bm{Y}} = \widetilde{\bm{A}}^{(:,\Theta)}\bm{\Omega}^* =
\left[\begin{array}{c}
\bm{\Sigma}_1\bm{\Omega}_1 \\
\bm{\Sigma}_2\bm{\Omega}_2
\end{array}\right].
$$
Denote $\bm{P}_{\bm{Y}}$ as the projection to the column space of $\bm{Y}$. We have
\begin{align*}
\|(\bm{I}-\bm{P}_{\bm{Y}})\bm{A}\|_F=&\|(\bm{I}-\bm{Q}\bm{Q}^*)\bm{A}\|_F\\
=&\|\bm{U}^*(\bm{I}-\bm{Q}\bm{Q}^*)\bm{U}\widetilde{\bm{A}}\|_F\\
=&\|(\bm{I}-(\bm{U}^*\bm{Q})(\bm{U}^*\bm{Q})^*)\widetilde{\bm{A}}\|_F
=\|(\bm{I}-\bm{P}_{\widetilde{\bm{Y}}})\widetilde{\bm{A}}\|_F.
\end{align*}
It suffices to prove that
$$
\|(\bm{I}-\bm{P}_{\widetilde{\bm{Y}}})\widetilde{\bm{A}}\|_F^2\leq \|\bm{\Sigma}_2\|_F^2+\|\bm{\Sigma}_2\bm{\Omega}_2\bm{\Omega}_1^{\dagger}\|_F^2.
$$

Second, we assume that $\bm{\Sigma}_2$ is not zero matrix. Otherwise,
$$
\text{range}(\widetilde{\bm{A}})
=\text{range}\left(\left[\begin{array}{c}
\bm{\Sigma}_1\bm{V}_1^* \\
\bm{0}
\end{array}\right]\right)
=\text{range}\left(\left[\begin{array}{c}
\bm{\Sigma}_1\bm{\Omega}_1 \\
\bm{0}
\end{array}\right]\right)
=\text{range}(\widetilde{\bm{Y}}),
$$
where the second equality holds since both $\bm{V}_1^*$ and $\bm{\Omega}_1$ have full row rank. As the result,
$$
\|(\bm{I}-\bm{P}_{\widetilde{\bm{Y}}})\widetilde{\bm{A}}\|_F=0,
$$
and the conclusion follows.

Next, we flatten out the top block of $\widetilde{\bm{Y}}$ to obtain
$$
\bm{Z} = \widetilde{\bm{Y}}\cdot\bm{\Omega}_1^{\dagger}\bm{\Sigma}_1^{-1}
=\left[\begin{array}{c}
\bm{I} \\
\bm{F}
\end{array}\right],\quad
\text{where} \quad
\bm{F} = \bm{\Sigma}_2\bm{\Omega}_2\bm{\Omega}_1^{\dagger}\bm{\Sigma}_1^{-1}.
$$
Since $\text{range}(\bm{Z})\subset\text{range}(\widetilde{\bm{Y}})$, 
$$
\|(\bm{I}-\bm{P}_{\widetilde{\bm{Y}}})\widetilde{\bm{A}}\|_F
\leq \|(\bm{I}-\bm{P}_{\bm{Z}})\widetilde{\bm{A}}\|_F
$$
by applying \cite[Proposition~8.5]{halko2011finding}. Taking squares,
\begin{align*}
\|(\bm{I}-\bm{P}_{\widetilde{\bm{Y}}})\widetilde{\bm{A}}\|_F^2
\leq\|(\bm{I}-\bm{P}_{\bm{Z}})\widetilde{\bm{A}}\|_F^2
= \text{Tr}(\widetilde{\bm{A}}^*(\bm{I}-\bm{P}_{\bm{Z}})\widetilde{\bm{A}})
= \text{Tr}(\bm{\Sigma}(\bm{I}-\bm{P}_{\bm{Z}})\bm{\Sigma}).
\end{align*}

Note that $\bm{Z}$ has full column rank,
$$
\bm{P}_{\bm{Z}} = \bm{Z}(\bm{Z}^*\bm{Z})^{-1}\bm{Z}^*=
\left[\begin{array}{c}
\bm{I} \\
\bm{F}
\end{array}\right]
(\bm{I}+\bm{F}^*\bm{F})^{-1}
\left[\begin{array}{c}
\bm{I} \\
\bm{F}
\end{array}\right]^*,
$$
and $\bm{I}-\bm{P}_{\bm{Z}}$ is equal to
$$
\left[\begin{array}{cc}
\bm{I} - (\bm{I}+\bm{F}^*\bm{F})^{-1} & -(\bm{I}+\bm{F}^*\bm{F})^{-1}\bm{F}^* \\
-\bm{F}(\bm{I}+\bm{F}^*\bm{F})^{-1} & \bm{I}-\bm{F}(\bm{I}+\bm{F}^*\bm{F})^{-1}\bm{F}^*
\end{array}\right].
$$
The top left block satisfies
$$
\bm{I} - (\bm{I}+\bm{F}^*\bm{F})^{-1} \preceq \bm{F}^*\bm{F}
$$
by applying \cite[Proposition~8.2]{halko2011finding}, and the bottom right block satisfies
$$
\bm{I}-\bm{F}(\bm{I}+\bm{F}^*\bm{F})^{-1}\bm{F}^* \preceq \bm{I}
$$
since $\bm{F}(\bm{I}+\bm{F}^*\bm{F})^{-1}\bm{F}^*\succeq \bm{0}$. Denote $\bm{B}=-(\bm{I}+\bm{F}^*\bm{F})^{-1}\bm{F}^*$,
$$
\bm{I}-\bm{P}_{\bm{Z}}\preceq 
\left[\begin{array}{cc}
\bm{F}^*\bm{F} & \bm{B} \\
\bm{B}^* & \bm{I}
\end{array}\right],
$$
and consequently,
$$
\bm{\Sigma}(\bm{I}-\bm{P}_{\bm{Z}})\bm{\Sigma}\preceq
\left[\begin{array}{cc}
\bm{\Sigma}_1\bm{F}^*\bm{F}\bm{\Sigma}_1 & \bm{\Sigma}_1\bm{B}\bm{\Sigma}_2 \\
\bm{\Sigma}_2\bm{B}^*\bm{\Sigma}_1 & \bm{\Sigma}_2^2
\end{array}\right].
$$
The last step is to note that
\begin{align*}
\text{Tr}(\bm{\Sigma}(\bm{I}-\bm{P}_{\bm{Z}})\bm{\Sigma})\leq&\text{Tr}\left(
\left[\begin{array}{cc}
\bm{\Sigma}_1\bm{F}^*\bm{F}\bm{\Sigma}_1 & \bm{\Sigma}_1\bm{B}\bm{\Sigma}_2 \\
\bm{\Sigma}_2\bm{B}^*\bm{\Sigma}_1 & \bm{\Sigma}_2^2
\end{array}\right]
\right)\\
=&\text{Tr}(\bm{\Sigma}_1\bm{F}^*\bm{F}\bm{\Sigma}_1)+\text{Tr}(\bm{\Sigma}_2^2)\\
=&\|\bm{F}\bm{\Sigma}_1\|_F^2+\|\bm{\Sigma}_2\|_F^2\\
=&\|\bm{\Sigma}_2\bm{\Omega}_2\bm{\Omega}_1^{\dagger}\|_F^2+\|\bm{\Sigma}_2\|_F^2.
\end{align*} 
\section{Proof of Theorem 2}

The final approximation $[[\hat{\bm{A}}]]_r$ is the best rank $r$ approximation of $\hat{\bm{A}}$. Note that
\begin{align*}
    \|\bm{A}-[[\hat{\bm{A}}]]_r\|_F 
    \leq& \|\bm{A}-\hat{\bm{A}}\|_F+\|\hat{\bm{A}}-[[\hat{\bm{A}}]]_r\|_F \\
    \leq& \|\bm{A}-\hat{\bm{A}}\|_F+\|\hat{\bm{A}}-[[\bm{A}]]_r\|_F \\
    \leq& 2\|\bm{A}-\hat{\bm{A}}\|_F+\|\bm{A}-[[\bm{A}]]_r\|_F\\
    =&2\|\bm{A}-\hat{\bm{A}}\|_F+\|\bm{\Sigma}_2\|_F,
\end{align*}
where the second inequality holds since $[[\hat{\bm{A}}]]_r$ is the best rank $r$ approximation to $\hat{\bm{A}}$. We could bound $\|\bm{A}-[[\hat{\bm{A}}]]_r\|_F$ by bounding the initial approximation error $\|\bm{A}-\hat{\bm{A}}\|_F$.
\begin{align*}
    \|\bm{A}-\hat{\bm{A}}\|_F^2
    =& \|\bm{A}-\bm{Q}\bm{C}\bm{P}^*\|_F^2 \\
    =& \|\bm{A}-\bm{Q}\bm{Q}^*\bm{A}\bm{P}\bm{P}^*+\bm{Q}(\bm{Q}^*\bm{A}\bm{P}-\bm{C})\bm{P}^*\|_F^2\\
    =& \|\bm{A}-\bm{Q}\bm{Q}^*\bm{A}\bm{P}\bm{P}^*\|_F^2+\|\bm{Q}(\bm{C}-\bm{Q}^*\bm{A}\bm{P})\bm{P}^*\|_F^2\\
    =&\underbrace{\|\bm{A}-\bm{Q}\bm{Q}^*\bm{A}\bm{P}\bm{P}^*\|_F^2}_{T_2}+\underbrace{\|\bm{C}-\bm{Q}^*\bm{A}\bm{P}\|_F^2}_{T_3},
\end{align*}
where the third equality is due to
$$
\langle\bm{A}-\bm{Q}\bm{Q}^*\bm{A}\bm{P}\bm{P}^*,\bm{Q}(\bm{C}-\bm{Q}^*\bm{A}\bm{P})\bm{P}^*\rangle=\langle\bm{Q}^*(\bm{A}-\bm{Q}\bm{Q}^*\bm{A}\bm{P}\bm{P}^*)\bm{P},\bm{C}-\bm{Q}^*\bm{A}\bm{P}\rangle=0.
$$

\textbf{Bound for $T_2$.} Note that
\begin{align*}
T_2
=&\|\bm{A}(\bm{I}-\bm{P}\bm{P}^*)+(\bm{I}-\bm{Q}\bm{Q}^*)\bm{A}\bm{P}\bm{P}^*\|_F^2\\
=&\|\bm{A}(\bm{I}-\bm{P}\bm{P}^*)\|_F^2+\|(\bm{I}-\bm{Q}\bm{Q}^*)\bm{A}\bm{P}\bm{P}^*\|_F^2\\
=&\|\bm{A}(\bm{I}-\bm{P}\bm{P}^*)\|_F^2+\|(\bm{I}-\bm{Q}\bm{Q}^*)\bm{A}\bm{P}\|_F^2\\
=&\|(\bm{U}_2\bm{\Sigma}_2\bm{V}_2^*)(\bm{I}-\bm{P}\bm{P}^*)\|_F^2+\|(\bm{U}_1\bm{\Sigma}_1\bm{V}_1^*)(\bm{I}-\bm{P}\bm{P}^*)\|_F^2+\|(\bm{I}-\bm{Q}\bm{Q}^*)\bm{A}\bm{P}\|_F^2\\
\leq&\|\bm{\Sigma}_2\|_F^2+\|(\bm{U}_1\bm{\Sigma}_1\bm{V}_1^*)(\bm{I}-\bm{P}\bm{P}^*)\|_F^2+\|(\bm{I}-\bm{Q}\bm{Q}^*)\bm{A}\bm{P}\|_F^2,
\end{align*}
where the second equality is due to $$\langle\bm{A}(\bm{I}-\bm{P}\bm{P}^*),(\bm{I}-\bm{Q}\bm{Q}^*)\bm{A}\bm{P}\bm{P}^*\rangle=0,$$ 
and the fourth equality is due to $$\langle(\bm{U}_2\bm{\Sigma}_2\bm{V}_2^*)(\bm{I}-\bm{P}\bm{P}^*),(\bm{U}_1\bm{\Sigma}_1\bm{V}_1^*)(\bm{I}-\bm{P}\bm{P}^*)\rangle=0.$$ 
Denote $\bm{\Gamma}_1=\bm{\Gamma}\bm{U}_1^{(\Delta,:)}$, and $\bm{\Gamma}_2=\bm{\Gamma}\bm{U}_2^{(\Delta,:)}$. Following a similar argument as the proof of Lemma~\ref{lem:deterministic},
$$
\|(\bm{U}_1\bm{\Sigma}_1\bm{V}_1^*)(\bm{I}-\bm{P}\bm{P}^*)\|_F^2\leq\|\bm{\Gamma}_1^{\dagger}\bm{\Gamma}_2\bm{\Sigma}_2\|_F^2.$$
Therefore, we arrive at
\begin{align*}
T_2
\leq& \|\bm{\Sigma}_2\|_F^2+\|\bm{\Gamma}_1^{\dagger}\bm{\Gamma}_2\bm{\Sigma}_2\|_F^2+\|(\bm{I}-\bm{Q}\bm{Q}^*)\bm{A}\bm{P}\|_F^2\\
\leq&\|\bm{\Sigma}_2\|_F^2+\|\bm{\Gamma}_1^{\dagger}\bm{\Gamma}_2\bm{\Sigma}_2\|_F^2+\|(\bm{I}-\bm{Q}\bm{Q}^*)\bm{A}\|_F^2\\
\leq&\|\bm{\Sigma}_2\|_F^2+\|\bm{\Gamma}_1^{\dagger}\bm{\Gamma}_2\bm{\Sigma}_2\|_F^2+\|\bm{\Sigma}_2\|_F^2+\|\bm{\Sigma}_2\bm{\Omega}_2\bm{\Omega}_1^{\dagger}\|_F^2\\
=& 2\|\bm{\Sigma}_2\|_F^2+\|\bm{\Gamma}_1^{\dagger}\bm{\Gamma}_2\bm{\Sigma}_2\|_F^2+\|\bm{\Sigma}_2\bm{\Omega}_2\bm{\Omega}_1^{\dagger}\|_F^2.
\end{align*}

\textbf{Bound for $T_3$.} 
We  decompose $\bm{C}-\bm{Q}^*\bm{A}\bm{P}$ into several parts, and bound them separately. Denote the complement of $\bm{Q}$ and $\bm{P}$ to be $\bm{Q}_{\perp}$ and $\bm{P}_{\perp}$, respectively. Define
$$
\bm{\Phi}_1 = \bm{\Phi}\bm{Q}^{(\Delta^{\prime},:)}\quad\text{and}\quad
\bm{\Phi}_2 = \bm{\Phi}\bm{Q}_{\perp}^{(\Delta^{\prime},:)},
$$
$$
\bm{\Psi}_1 = \bm{\Psi}\bm{P}^{(\Theta^{\prime},:)}\quad\text{and}\quad
\bm{\Psi}_2 = \bm{\Psi}\bm{P}_{\perp}^{(\Theta^{\prime},:)}.
$$
\begin{lemma}\label{lem:decomposition}
Assume that $\bm{\Phi}_1$ and $\bm{\Psi}_1$ are of full column rank. Then
\begin{align*}
\bm{C}-\bm{Q}^*\bm{A}\bm{P}=\bm{\Phi}_1^{\dagger}\bm{\Phi}_2(\bm{Q}_{\perp}^*\bm{A}\bm{P})+(\bm{Q}^*\bm{A}\bm{P}_{\perp})\bm{\Psi}_2^*(\bm{\Psi}_1^{\dagger})^*+\bm{\Phi}_1^{\dagger}\bm{\Phi}_2(\bm{Q}_{\perp}^*\bm{A}\bm{P}_{\perp})\bm{\Psi}_2^*(\bm{\Psi}_1^{\dagger})^*.
\end{align*}
\end{lemma}
The proof of Lemma~\ref{lem:decomposition} is the same as \cite[Lemma~A.3]{tropp2019streaming}. For completeness, we include its proof in Section~\ref{sec:Lemma7}.

\begin{remark*}
To ensure that $\bm{\Phi}_1$ and $\bm{\Psi}_1$ are of full column rank, one way is to argue that $\bm{Q}$ and $\bm{P}$ are also incoherent, i.e.,
\begin{align}\label{crux}
\max_i\left\{\|\bm{Q}^{(i,:)}\|_2\right\} \leq \sqrt{\frac{\mu^{\prime} k}{M}}
\quad\text{and}\quad
\max_j\left\{\|\bm{P}^{(j,:)}\|_2\right\} \leq \sqrt{\frac{\nu^{\prime} k}{N}}
\end{align}
for some $\mu^{\prime}\ll M$, $\nu^{\prime}\ll N$, and apply Lemma~\ref{lem:incoherence}. As shown in the following table, we empirically observe that \eqref{crux} holds, and the incoherence parameters $\mu^{\prime}=O(\mu)$, $\nu^{\prime}=O(\nu)$. Table $6$ also confirms that $(\mu,\nu)$ and $(\mu^{\prime},\nu^{\prime})$ are indeed small, compared to matrix size $(M,N)$, for all the datasets.

\begin{table}[H]
\begin{center}
\small
\begin{tabular}{ccccc}
\hline
\multicolumn{1}{|c|}{} & \multicolumn{1}{c|}{Yale Face} & \multicolumn{1}{c|}{Cardiac MRI} & \multicolumn{1}{c|}{BR1003} & \multicolumn{1}{c|}{Video} \\
\hline
\multicolumn{1}{|c|}{$\mu$} & \multicolumn{1}{c|}{4.1137} & \multicolumn{1}{c|}{127.5935} & \multicolumn{1}{c|}{5.4270} & \multicolumn{1}{c|}{20.3505} \\
\hline
\multicolumn{1}{|c|}{$\mu^{\prime}$} & \multicolumn{1}{c|}{5.9454} & \multicolumn{1}{c|}{159.8982} & \multicolumn{1}{c|}{5.2750} & \multicolumn{1}{c|}{22.8109} \\
\hline
\multicolumn{1}{|c|}{$\nu$} & \multicolumn{1}{c|}{2.7068} & \multicolumn{1}{c|}{2.1507} & \multicolumn{1}{c|}{32.8387} & \multicolumn{1}{c|}{14.0194}  \\
\hline
\multicolumn{1}{|c|}{$\nu^{\prime}$} & \multicolumn{1}{c|}{4.1355} & \multicolumn{1}{c|}{2.2260} & \multicolumn{1}{c|}{57.4357} & \multicolumn{1}{c|}{6.7429}  \\
\hline
\end{tabular}
\caption{Comparison of $(\mu,\nu)$ and $(\mu^{\prime},\nu^{\prime})$ on different datasets. The parameters $(r,k,s,p)$ for each dataset are chosen as the same values as used to show all the visual comparisons.}
\end{center}
\end{table}
  However, we find that our theoretical bounds of $(\mu^{\prime},\nu^{\prime})$ in \eqref{crux} are currently difficult to prove. 
  One way to get around this is to slightly modify the algorithm: after building the left and right sketches, we calculate the QR factorizations, and sample $m^{\prime}$ row indices and $n^{\prime}$ column indices based on the actual incoherence parameters, estimated from the row norms, of $\bm{Q}$ and $\bm{P}$, respectively. Then we sample $m^{\prime}\approx O(\mu^{\prime}k\log k)$ row indices, and $n^{\prime}\approx O(\nu^{\prime}k\log k)$ column indices to get the data in the intersection.
\end{remark*}

\begin{align*}
    T_3
    =&\|\bm{\Phi}_1^{\dagger}\bm{\Phi}_2(\bm{Q}_{\perp}^*\bm{A}\bm{P})+
    (\bm{Q}^*\bm{A}\bm{P}_{\perp})\bm{\Psi}_2^*(\bm{\Psi}_1^{\dagger})^*+\bm{\Phi}_1^{\dagger}\bm{\Phi}_2(\bm{Q}_{\perp}^*\bm{A}\bm{P}_{\perp})\bm{\Psi}_2^*(\bm{\Psi}_1^{\dagger})^*\|_F^2\\
    \leq&3\cdot(\underbrace{\|\bm{\Phi}_1^{\dagger}\bm{\Phi}_2(\bm{Q}_{\perp}^*\bm{A}\bm{P})\|_F^2}_{T_{31}}+
    \underbrace{\|(\bm{Q}^*\bm{A}\bm{P}_{\perp})\bm{\Psi}_2^*(\bm{\Psi}_1^{\dagger})^*\|_F^2}_{T_{32}}+\underbrace{\|\bm{\Phi}_1^{\dagger}\bm{\Phi}_2(\bm{Q}_{\perp}^*\bm{A}\bm{P}_{\perp})\bm{\Psi}_2^*(\bm{\Psi}_1^{\dagger})^*\|_F^2}_{T_{33}}).
\end{align*}
We bound the three terms separately. The bounds all follow from a similar argument as deriving the bound for $\|\bm{\Sigma}_2\bm{\Omega}_2\bm{\Omega}_1^{\dagger}\|_F^2$ in the proof of Theorem 1. Suppose $s\geq k+4$.
\begin{align*}
T_{31}\leq&\frac{6e^2}{q}\cdot\frac{s^{1+6/(s-k+1)}}{(s-k+1)^2}\cdot\left(\sqrt{s}\|\bm{Q}_{\perp}^*\bm{A}\bm{P}\|_F+\sqrt{6\log s}\|\bm{Q}_{\perp}^*\bm{A}\bm{P}\|_2\right)^2\\
\leq&\frac{6e^2}{q}\cdot\frac{s^{1+6/(s-k+1)}}{(s-k+1)^2}\cdot\left(\sqrt{s}+\sqrt{6\log s}\right)^2\cdot\|\bm{Q}_{\perp}^*\bm{A}\bm{P}\|_F^2\\
=& C(q,s,k)\cdot\|\bm{Q}_{\perp}^*\bm{A}\bm{P}\|_F^2
\end{align*}
with probability at least $1-\frac{2}{s^3}$, where
$$
C(q,s,k)=\frac{6e^2}{q}\cdot\frac{s^{1+6/(s-k+1)}}{(s-k+1)^2}\cdot\left(\sqrt{s}+\sqrt{6\log s}\right)^2.
$$
Similarly, with probability at least $1-\frac{2}{s^3}$,
$$
T_{32}\leq C(q,s,k)\cdot\|\bm{Q}^*\bm{A}\bm{P}_{\perp}\|_F^2.
$$
Apply the argument twice,
$$
T_{33} \leq C(q,s,k)\cdot\|(\bm{Q}_{\perp}^*\bm{A}\bm{P}_{\perp})\bm{\Psi}_2^*(\bm{\Psi}_1^{\dagger})^*\|_F^2
\leq (C(q,s,k))^2\cdot\|\bm{Q}_{\perp}^*\bm{A}\bm{P}_{\perp}\|_F^2
$$
with probability at least $1-\frac{2}{s^2}$. Combining the three estimates,
\begin{align*}
\frac{T_3}{3}
\leq & C(q,s,k)\cdot\left(\|\bm{Q}_{\perp}^*\bm{A}\bm{P}\|_F^2+\|\bm{Q}^*\bm{A}\bm{P}_{\perp}\|_F^2+C(q,s,k)\cdot\|\bm{Q}_{\perp}^*\bm{A}\bm{P}_{\perp}\|_F^2\right)\\
= &C(q,s,k)\cdot\left(\|\bm{Q}_{\perp}^*\bm{A}\bm{P}\|_F^2+\|\bm{Q}^*\bm{A}\bm{P}_{\perp}\|_F^2+\|\bm{Q}_{\perp}^*\bm{A}\bm{P}_{\perp}\|_F^2+(C(q,s,k)-1)\cdot\|\bm{Q}_{\perp}^*\bm{A}\bm{P}_{\perp}\|_F^2\right)\\
= & C(q,s,k)\cdot\left(\|\bm{A}-\bm{Q}\bm{Q}^*\bm{A}\bm{P}\bm{P}^*\|_F^2+(C(q,s,k)-1)\cdot\|\bm{Q}_{\perp}^*\bm{A}\bm{P}_{\perp}\|_F^2\right)\\
= & C(q,s,k)\cdot\|\bm{A}-\bm{Q}\bm{Q}^*\bm{A}\bm{P}\bm{P}^*\|_F^2+C_(q,s,k)(C(q,s,k)-1)\cdot\|\bm{Q}_{\perp}\bm{Q}_{\perp}^*\bm{A}\bm{P}_{\perp}\|_F^2\\
\leq & C(q,s,k)\cdot\|\bm{A}-\bm{Q}\bm{Q}^*\bm{A}\bm{P}\bm{P}^*\|_F^2+C(q,s,k)(C(q,s,k)-1)\cdot\|(\bm{I}-\bm{Q}\bm{Q}^*)\bm{A}\|_F^2.
\end{align*}
Thus for the square of the initial approximation error,
\begin{align*}
    \|\bm{A}-\hat{\bm{A}}\|_F^2
    \leq&(3C(q,s,k)+1)\cdot\|\bm{A}-\bm{Q}\bm{Q}^*\bm{A}\bm{P}\bm{P}^*\|_F^2\\
    &+3C(q,s,k)(C(q,s,k)-1)\cdot\|(\bm{I}-\bm{Q}\bm{Q}^*)\bm{A}\bm{P}\|_F^2\\
    \leq&(3C(q,s,k)+1)\cdot(2\|\bm{\Sigma}_2\|_F^2+\|\bm{\Gamma}_1^{\dagger}\bm{\Gamma}_2\bm{\Sigma}_2\|_F^2+\|\bm{\Sigma}_2\bm{\Omega}_2\bm{\Omega}_1^{\dagger}\|_F^2)\\
    &+3C(q,s,k)(C(q,s,k)-1)\cdot(\|\bm{\Sigma}_2\|_F^2+\|\bm{\Sigma}_2\bm{\Omega}_2\bm{\Omega}_1^{\dagger}\|_F^2).
\end{align*}
Taking the square root,
\begin{align*}
    \|\bm{A}-\hat{\bm{A}}\|_F
    \leq&\sqrt{3(C(q,s,k))^2+3C(q,s,k)+2}\cdot(C_1(p,k,r)\|\bm{\Sigma}_2\|_F+C_2(p,k,r)\|\bm{\Sigma}_2\|_2)\\
    \leq&(\sqrt{3}C(q,s,k)+\sqrt{2})\cdot(C_1(p,k,r)\|\bm{\Sigma}_2\|_F+C_2(p,k,r)\|\bm{\Sigma}_2\|_2)\\
    =&C_3(p,q,s,k,r)\cdot\|\bm{\Sigma}_2\|_F+C_4(p,q,s,k,r)\cdot\|\bm{\Sigma}_2\|_2,
\end{align*}
where
\begin{align*}
C_3(p,q,s,k,r)=&C_1(p,k,r)\cdot(\sqrt{3}C(q,s,k)+\sqrt{2}),\\
C_4(p,q,s,k,r)=&C_2(p,k,r)\cdot(\sqrt{3}C(q,s,k)+\sqrt{2}). 
\end{align*}

\subsection{Proof of Lemma~\ref{lem:decomposition}}
\label{sec:Lemma7}

Define $\bm{S}_{\Delta^{\prime}}\in\mathbb{R}^{m^{\prime}\times M}$ such that when left multiplied to a matrix, the result is equal to the rows with indices $\Delta^{\prime}$ of that matrix. Define $\bm{S}_{\Theta^{\prime}}$ similarly. The core sketch can be written as
\begin{align*}
\bm{Z}=&\bm{\Phi}\bm{A}^{(\Delta^{\prime},\Theta^{\prime})}\bm{\Psi}^*
=\bm{\Phi}\bm{S}_{\Delta^{\prime}}\bm{A}\bm{S}_{\Theta^{\prime}}^*\bm{\Psi}^*\\
=&\bm{\Phi}\bm{S}_{\Delta^{\prime}}(\bm{A}-\bm{Q}\bm{Q}^*\bm{A}\bm{P}\bm{P}^*)\bm{S}_{\Theta^{\prime}}^*\bm{\Psi}^*+(\bm{\Phi}\bm{S}_{\Delta^{\prime}}\bm{Q})\bm{Q}^*\bm{A}\bm{P}(\bm{P}^*\bm{S}_{\Theta^{\prime}}^*\bm{\Psi}^*)\\
=&\bm{\Phi}\bm{S}_{\Delta^{\prime}}(\bm{A}-\bm{Q}\bm{Q}^*\bm{A}\bm{P}\bm{P}^*)\bm{S}_{\Theta^{\prime}}^*\bm{\Psi}^*+(\bm{\Phi}\bm{Q}^{(\Delta^{\prime},:)})\bm{Q}^*\bm{A}\bm{P}(\bm{\Psi}\bm{P}^{(\Theta^{\prime},:)})^*.
\end{align*}
Left multiply by $\bm{\Phi}_1^{\dagger}$ and right-multiply by $(\bm{\Psi}_1^{\dagger})^*$,
\begin{align*}
\bm{C}=&\bm{\Phi}_1^{\dagger}\bm{\Phi}\bm{S}_{\Delta^{\prime}}(\bm{A}-\bm{Q}\bm{Q}^*\bm{A}\bm{P}\bm{P}^*)\bm{S}_{\Theta^{\prime}}^*\bm{\Psi}^*(\bm{\Psi}_1^{\dagger})^*+\bm{Q}^*\bm{A}\bm{P}.
\end{align*}
Notice that
$$
\bm{\Phi}_1^{\dagger}\bm{\Phi}\bm{S}_{\Delta^{\prime}}
=\bm{\Phi}_1^{\dagger}\bm{\Phi}\bm{S}_{\Delta^{\prime}}\bm{Q}\bm{Q}^*+\bm{\Phi}_1^{\dagger}\bm{\Phi}\bm{S}_{\Delta^{\prime}}\bm{Q}_{\perp}\bm{Q}_{\perp}^*=\bm{Q}^*+\bm{\Phi}_1^{\dagger}\bm{\Phi}_2\bm{Q}_{\perp}^*,
$$
$$
\bm{S}_{\Theta^{\prime}}^*\bm{\Psi}^*(\bm{\Psi}_1^{\dagger})^*
=\bm{P}\bm{P}^*\bm{S}_{\Theta^{\prime}}^*\bm{\Psi}^*(\bm{\Psi}_1^{\dagger})^*+
\bm{P}_{\perp}\bm{P}_{\perp}^*\bm{S}_{\Theta^{\prime}}^*\bm{\Psi}^*(\bm{\Psi}_1^{\dagger})^*
=\bm{P}+\bm{P}_{\perp}\bm{\Psi}_2^*(\bm{\Psi}_1^{\dagger})^*.
$$
Combining all the pieces,
\begin{align*}
\bm{C}-\bm{Q}^*\bm{A}\bm{P}=&\left(\bm{Q}^*+\bm{\Phi}_1^{\dagger}\bm{\Phi}_2\bm{Q}_{\perp}^*\right)(\bm{A}-\bm{Q}\bm{Q}^*\bm{A}\bm{P}\bm{P}^*)\left(\bm{P}+\bm{P}_{\perp}\bm{\Psi}_2^*(\bm{\Psi}_1^{\dagger})^*\right)\\
=&\bm{\Phi}_1^{\dagger}\bm{\Phi}_2(\bm{Q}_{\perp}^*\bm{A}\bm{P})+(\bm{Q}^*\bm{A}\bm{P}_{\perp})\bm{\Psi}_2^*(\bm{\Psi}_1^{\dagger})^*+\bm{\Phi}_1^{\dagger}\bm{\Phi}_2(\bm{Q}_{\perp}^*\bm{A}\bm{P}_{\perp})\bm{\Psi}_2^*(\bm{\Psi}_1^{\dagger})^*.
\end{align*} 
\section{More Numerics}

We use a Navier Stokes simulated flow system which generates vorticity patterns for an incompressible fluid under certain initial and boundary conditions. For each point of the $100\times 50$ grids, fluid  velocity values in both $x$ and $y$ directions for $200$ time instances. As an example we demonstrate our results on the $y$ component of the fluid velocity, captured in the data matrix $\bm{A}\in\mathbb{R}^{5000\times 200}$. Based on the scree plot, we choose $r=7$. The optimal $err$ is 0.0016.    

\begin{figure}[H]
\label{fig:scree_NS}
    \centering
    \includegraphics[width=0.35\textwidth]{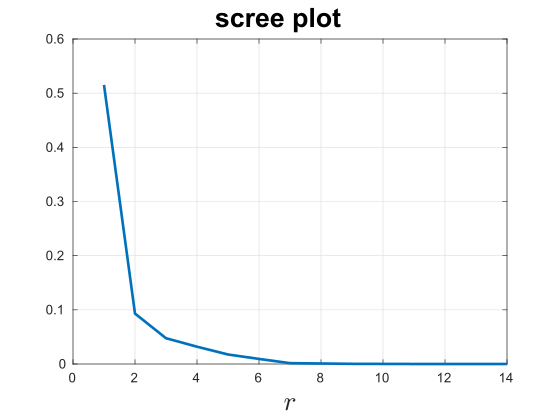}
    \caption{Scree plot for Navier Stokes dataset.}
\end{figure}

\begin{table}[H]
\begin{center}
\small
\begin{tabular}{ccc}
\hline
\multicolumn{1}{|c|}{} & \multicolumn{1}{c|}{SketchySVD} & \multicolumn{1}{c|}{SketchyCoreSVD} \\
\hline
\multicolumn{1}{|c|}{$p$} & \multicolumn{1}{c|}{-} & \multicolumn{1}{c|}{0.3} \\
\hline
\multicolumn{1}{|c|}{$err$} & \multicolumn{1}{c|}{0.0016} & \multicolumn{1}{c|}{0.0016} \\
\hline
\multicolumn{1}{|c|}{time (sec)} & \multicolumn{1}{c|}{0.0097} & \multicolumn{1}{c|}{0.0062} \\
\hline
\end{tabular}
\caption{Performance comparisons for Navier Stokes dataset.}
\end{center}
\end{table}

We choose $k=4r+1=29$ and $s=2k+1=59$ for both methods. In  Table $7$, we see that the error of SketchySVD is the same as the optimal error. SketchyCoreSVD can thus achieve the same error bound in less time. In the visual comparison Figure $10$, we show that the singular vector(s) can be estimated accurately, and with sampling ratio $p=30\%$ for SketchyCoreSVD.
Table $8$ further verifies that the PSNR ratios of the computed singular vectors are very high, and that SketchyCoreSVD achieves those ratios with only sampling ratio $p=30\%$.

\begin{table}[H]
\begin{center}
\small
\begin{tabular}{ccccccc}
\hline
\multicolumn{1}{|c|}{} & \multicolumn{1}{c|}{1} & \multicolumn{1}{c|}{2} & \multicolumn{1}{c|}{3} & \multicolumn{1}{c|}{4} & \multicolumn{1}{c|}{5} & \multicolumn{1}{c|}{6} \\
\hline
\multicolumn{1}{|c|}{SketchySVD} & \multicolumn{1}{c|}{115.3884} & \multicolumn{1}{c|}{116.2469} & \multicolumn{1}{c|}{112.3364} & \multicolumn{1}{c|}{96.8527} & \multicolumn{1}{c|}{97.7459} & \multicolumn{1}{c|}{88.7866} \\
\hline
\multicolumn{1}{|c|}{SketchyCoreSVD} & \multicolumn{1}{c|}{107.3097} & \multicolumn{1}{c|}{107.2117} & \multicolumn{1}{c|}{98.2154} & \multicolumn{1}{c|}{84.5132} & \multicolumn{1}{c|}{83.6819} & \multicolumn{1}{c|}{85.8143} \\
\hline
\end{tabular}
\caption{Comparison of PSNR ratios for the first six left singular vectors computed on the Navier Stokes dataset, where $p=0.3$ for our SketchyCoreSVD method.}
\end{center}
\end{table}

\begin{figure}[H]
    \centering
    \includegraphics[width=0.3\textwidth]{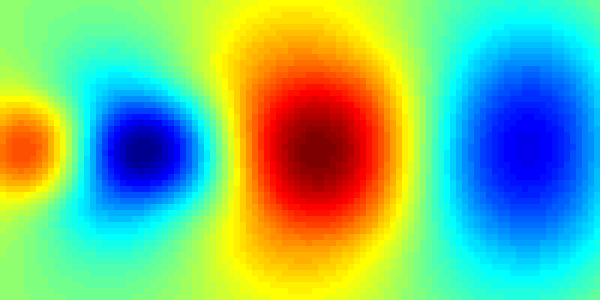}
    \includegraphics[width=0.3\textwidth]{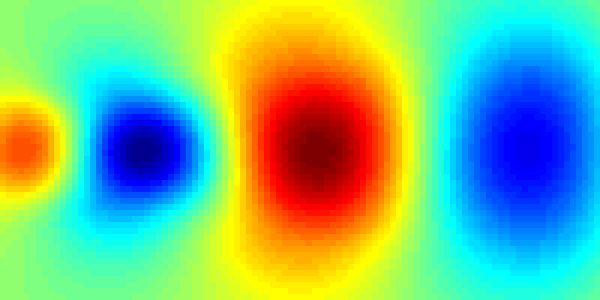}
    \includegraphics[width=0.3\textwidth]{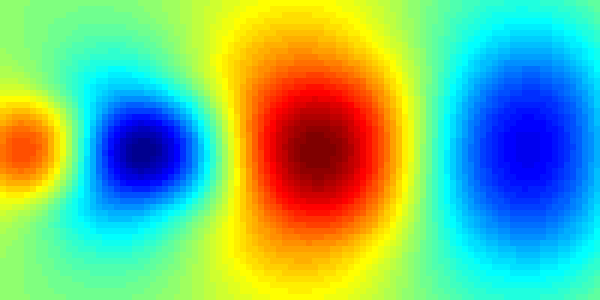}\\
    \includegraphics[width=0.3\textwidth]{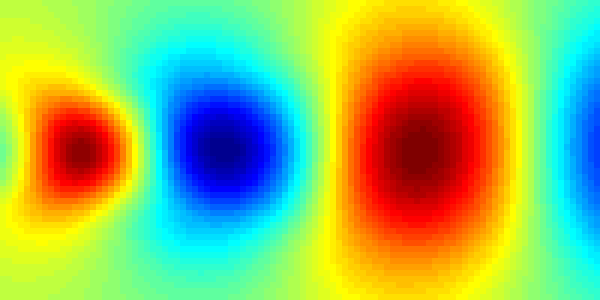}
    \includegraphics[width=0.3\textwidth]{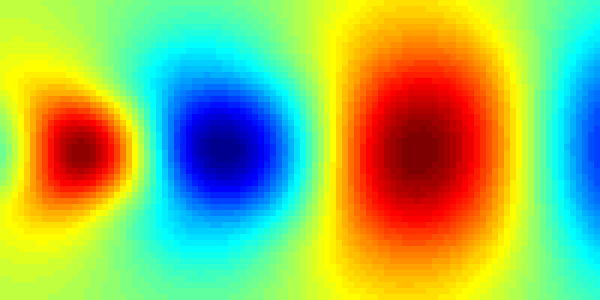}
    \includegraphics[width=0.3\textwidth]{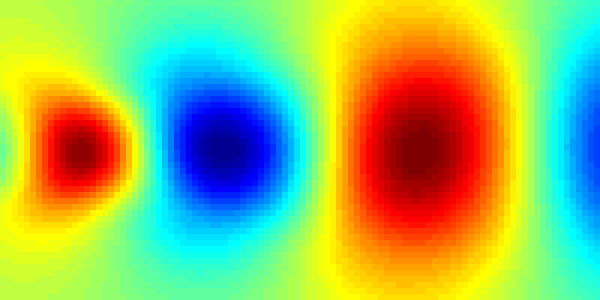}\\
    \includegraphics[width=0.3\textwidth]{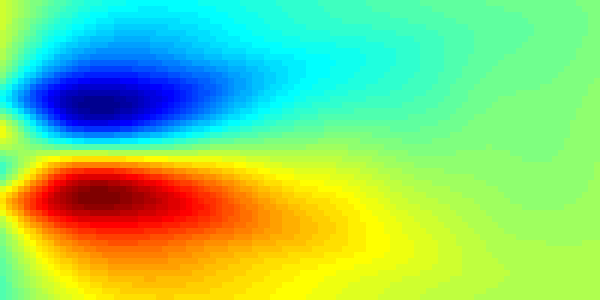}
    \includegraphics[width=0.3\textwidth]{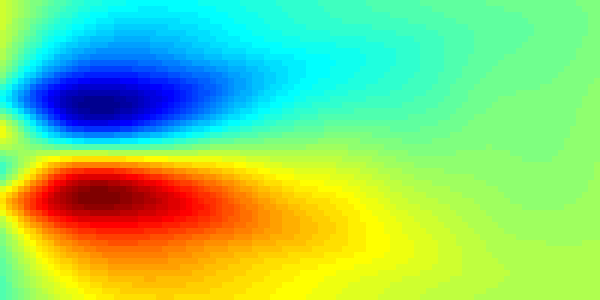}
    \includegraphics[width=0.3\textwidth]{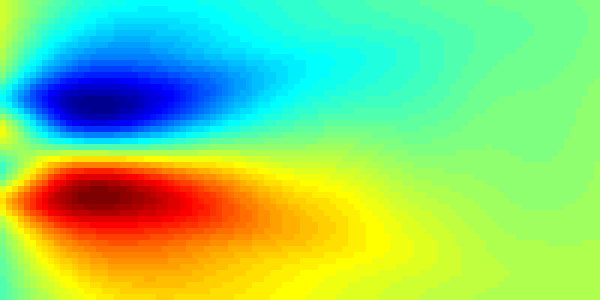}\\
    \includegraphics[width=0.3\textwidth]{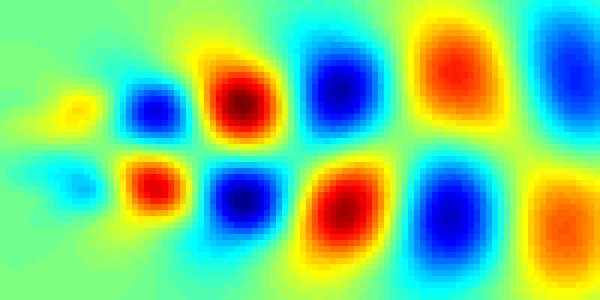}
    \includegraphics[width=0.3\textwidth]{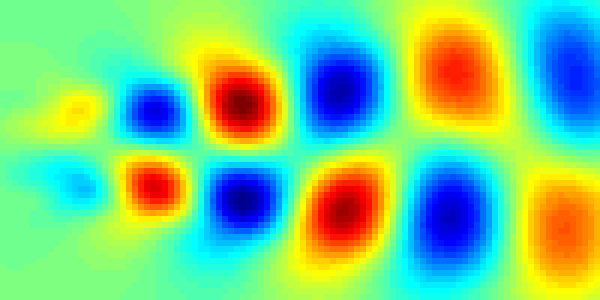}
    \includegraphics[width=0.3\textwidth]{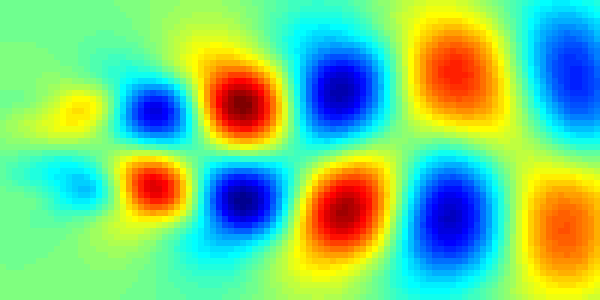}\\
    \includegraphics[width=0.3\textwidth]{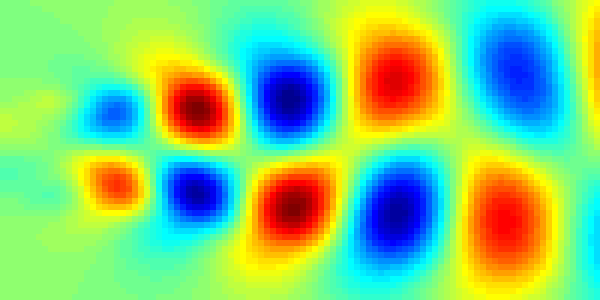}
    \includegraphics[width=0.3\textwidth]{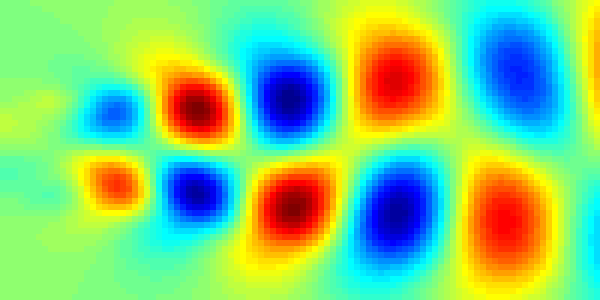}
    \includegraphics[width=0.3\textwidth]{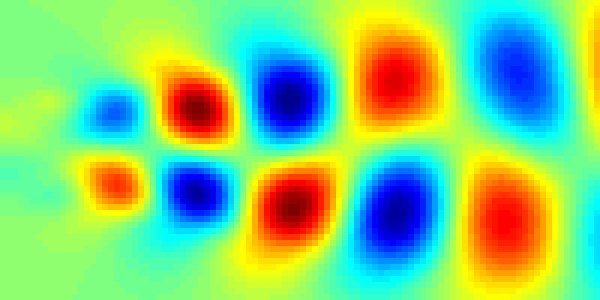}\\
    \includegraphics[width=0.3\textwidth]{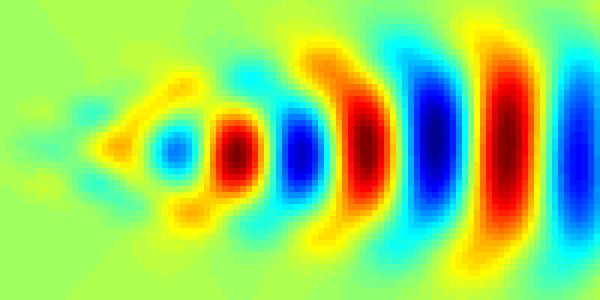}
    \includegraphics[width=0.3\textwidth]{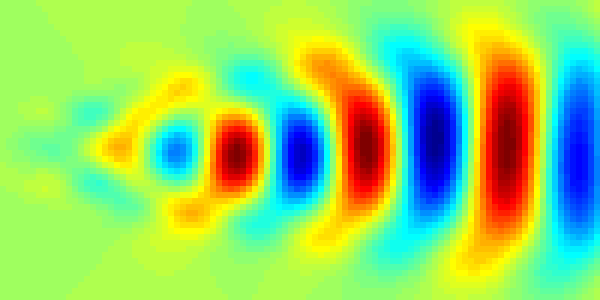}
    \includegraphics[width=0.3\textwidth]{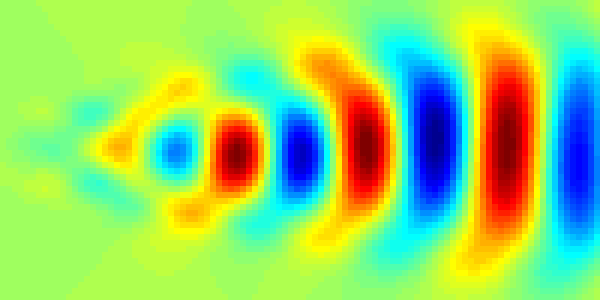}\\
    \caption{Visual comparison for the first six left singular vectors computed on the Navier Stokes dataset. Ground truth (left column), SketchySVD (middle column), SketchyCoreSVD (right column, $p=0.3$). For each subfigure, $x\in[1,8]$ (horizontal), and $y\in[-2,2]$ (vertical).}
\end{figure} 

\end{document}